\documentclass[a4paper, 11pt, twoside]{article}
\usepackage{amssymb, amsmath}
\usepackage{amsthm, amsfonts,mathrsfs}
\usepackage[T1]{fontenc}
\usepackage{times}
\usepackage{graphicx}
\usepackage{subfigure}
\usepackage{cite}
\usepackage{epsfig}
\usepackage{color}
\usepackage[all]{xypic}
\newtheorem{theorem}{Theorem}[section]

\newtheorem{remark}{Remark}[section]

\numberwithin{equation}{section}

\graphicspath{{image/}}

\begin{document} % (the body of the article begins)
%\begin{CJK*}{GBK}{song}
%  (title)
%\pagecolor{blue}
\title{{\large  \textbf{Global smooth solution for the 3D generalized tropical climate model with partial viscosity and damping}}\thanks{\footnotesize{
{The  authors are supported by the National Natural Science Foundation of China (No. 12271293, No. 11901342 and  No. 11701269),
Natural Science Foundation of Shandong Province (No. ZR2023MA002), Natural Science Foundation of Jiangsu Province (No. BK20231301), the project of Youth Innovation Team of Universities of Shandong Province (No. 2023KJ204).}
			}}}
\author{\small{Hui Liu$^1$\thanks{liuhuinanshi@qfnu.edu.cn}~~~~Chengfeng Sun$^{2}$\thanks{Corresponding author, sch200130@163.com}
~~~~Mei Li$^{2}$\thanks{limei@nufe.edu.cn}}\\
	\footnotesize{$1$. School of Mathematical Sciences, Qufu Normal University, Qufu, Shandong 273165, PR China;}
	\\
    \footnotesize{$2$. School of Applied Mathematics, Nanjing University of Finance and Economics, Nanjing, 210023,  PR China}
    }
\date{}
\maketitle
\noindent{\small{\hspace{1.1cm} }}
	%\begin{abstract}
	\\
\noindent \textbf{Abstract~~~}    The three-dimensional generalized tropical climate model with partial viscosity and damping is considered in this paper. Global well-posedness of solutions of the three-dimensional generalized tropical climate model with partial viscosity and damping is proved for $\alpha\geq\frac{3}{2}$ and $\beta\geq4$. Global smooth  solution of the three-dimensional generalized tropical climate model with partial viscosity and damping is proved in $H^s(\mathbb R^3)$ $(s>2)$  for $\alpha\geq\frac{3}{2}$ and $4\leq\beta\leq5$.
\\[2mm]
\textbf{Key words~~~}Tropical climate model; Damping; Global smooth solution
\\[2mm]
\\
\textbf{2020 Mathematics Subject Classification~~~}35B65, 35Q35, 35Q86, 76D03
%\end{abstract}
%\CenterWallPaper{1}{bg-a4-odd.jpg} %居中的背景图片，出现在文档的每一页上
\section {Introduction}
In this paper, we consider the following three-dimensional (3D) generalized tropical climate model with partial viscosity and damping:
 \begin{equation}\label{1}\begin{array}{l}
\left\{
\begin{array}{l}
\partial_{t}u+(u\cdot\nabla)u-\Delta_{h}u+\lvert u\rvert^{\beta-1}u+\nabla\cdot(v\otimes v)+\nabla p=0,\\
\partial_{t}v+(u\cdot\nabla)v+(-\Delta)^{\alpha}v+(v\cdot\nabla)u+\nabla\theta=0,\\
\partial_{t}\theta+(u\cdot\nabla)\theta-\Delta \theta+\nabla\cdot v=0,\\
\nabla\cdot u=0,\\
u(x,0)=u_{0}(x),v(x,0)=v_{0}(x),\theta(x,0)=\theta_{0}(x),
\end{array}
\right.
\end{array}\end{equation}
for $x\in\mathbb{R}^{3}$ and $t\geq0$. $u$ denotes the barotropic mode, $v$ denotes the first baroclinic mode of the velocity field, respectively. $p$ denotes the scalar pressure and $\theta$ denotes the scalar temperature, respectively. $\beta\geq1$ and $\alpha>0$ are real parameters. $v\otimes v$ is the standard tensor notation. $\Delta_{h}:=\partial_{1}^{2}+\partial_{2}^{2}$  and $\partial_{i}$ is the partial derivative in the direction $x_{i}$ and $i=1,2,3$. The fractional Laplacian operator $(-\Delta)^{\alpha}$ is defined via the Fourier transform
\begin{align*}
\widehat{(-\Delta)^{\alpha}}f(\xi)=|\xi|^{2\alpha}\widehat{f}(\xi).
\end{align*}

The tropical climate model was investigated by many authors in \cite{frierson,li,li1}.  Global well-posedness of strong solutions of the two-dimensional tropical climate model was proved by introducing the  pseudo-baroclinic velocity in \cite{li}. Meanwhile, the authors in \cite{li1}  proved the existence and uniqueness of global strong solutions of the two-dimensional tropical atmosphere model with moisture. By using the Littlewood-Paley theory and  a generalized commutator estimate, Wan proved the global small solutions of the two-dimensional tropical climate model without thermal diffusion in \cite{wan}. In \cite{ye}, Ye utilized the energy estimates  to prove the global regularity for a class of the two-dimensional tropical climate model. By constructing a new energy estimates in
Besov space, Ma and Wan in \cite{ma} proved the spectral analysis and  global well-posedness for a viscous two-dimensional tropical climate model with only a damp term. Global regularity of the  two-dimensional tropical climate model with fractional dissipation were proved in \cite{dong,dong2,dong1}. In \cite{yez}, Ye and Zhu proved the global strong solutions of the two-dimensional tropical climate model with temperature-dependent viscosity. Global smooth solution and asymptotic behavior of the two-dimensional temperature-dependent tropical climate model were proved in \cite{dong3,lixu}. By using the Littlewood-Paley decomposition and anisotropic inequalities, Paicu and Zhu in \cite{paicu} proved global regularity and large time behaviour of the two-dimensional  MHD and tropical climate model with horizontal dissipation in a strip domain $\mathbb{T}\times\mathbb{R}$.

The three-dimensional tropical climate model with damping were investigated in \cite{berti,yuan,yuan1}. Global strong solution of the three-dimensional tropical climate model with damping was proved in \cite{yuan,yuan1}. A regularity criterion of the three-dimensional tropical climate model with damping was obtained in \cite{berti}.

Recently, generalized magnetohydrodynamic (MHD) equations were extensively investigated in \cite{wu,wu1}. In \cite{wu1}, Wu proved the regularity criteria of the generalized MHD equations in Besov spaces. By introducing the corresponding vorticity field satisfies certain conditions, Zhou in \cite{zhou} proved the regularity criteria of the generalized viscous MHD equations. In addition,  Yamazaki in \cite{yamazaki} proved the regularity criteria of generalized MHD and Navier-Stokes equations. In \cite{zhao}, Zhao proved the higher-order derivative of solutions of the three-dimensional generalized Hall-MHD equations. By using the semigroup theory and  Hille-Yosida theory, Stefanov and  Hadadifard proved the sharp time decay rates of the two-dimensional generalized quasi-geostrophic equations and the Boussinesq equations in \cite{stefanov}.

Global regularity and time-periodic solutions of the two-dimensional incompressible magnetohydrodynamic  equations  with horizontal dissipation and horizontal magnetic diffusion were proved in \cite{cao,cao1,fan,sun}. Global existence of axisymmetric solutions of the three-dimensional Boussinesq equations with horizontal dissipation was proved in \cite{miao}. Global existence and regularity of the two-dimensional  Boussinesq equations with minimal dissipation and thermal diffusion were proved in \cite{wu2}. Cao, Li and Titi in \cite{cao2} proved the  local and the global well-posedness of strong solutions of the three-dimensional primitive equations with only horizontal dissipation in the periodic boundary conditions. Well-posedness of global solutions of the modified anisotropic three-dimensional Navier-Stokes equations were proved for $\beta>3$ in \cite{bessaih}. Existence and uniqueness of global solutions of the three-dimensional Boussinesq-MHD equations with partial viscosity and damping for $\beta\geq4$ with $\alpha=\frac{3}{2}$ in \cite{liu4}. Well-posedness and attractors of the fluids equation with damping were investigated in \cite{liu4,liu,liu1,liu2,liux}.

In order to get the existence and uniqueness of global solutions of the system (\ref{1}), we need to obtain the estimations $\|\partial_{3}u\|^{2}_{L^{2}}$ and
$\int_{\mathbb{R}^{3}}(u\cdot\nabla)\tilde{u}\cdot udx$ by using several anisotropic embedding and interpolation inequalities involving fractional derivatives.
If $\alpha=\frac32$, we will prove the Theorem \ref{the} by the similar method in \cite{liu4}. In \cite{liu4}, they haven't proved the global smooth solutions.   Due to loss the term $\nabla\cdot v=0$, we will prove the global smooth solutions of the system (\ref{1}) as the Theorem \ref{the1}. Compared with the classic Navier-Stokes equations, we replace the term $\Delta$ with the term $\Delta_{h}$. Therefore, in order to get the  global smooth solutions of the system (\ref{1}), we should overcome some difficulties and introduce several anisotropic embedding, such as the estimations $\|\partial_{t}u\|^{2}_{L^{2}}$ and $\|\partial_{t}v\|^{2}_{L^{2}}$ and $\|u\|_{L^{3(\beta+1)}}^{\beta+1}$ and $\|u\|_{L^{\infty}}^{\beta+1}$ and $\|\Lambda^{s}u\|^{2}_{L^{2}}$ and $\|\Lambda^{s}v\|^{2}_{L^{2}}$ and $\|\Lambda^{s}\theta\|^{2}_{L^{2}}$ $(s>2)$, and so on.

This paper is organized as follows. In section \ref{s1}, we give the basic definitions and main results. In section \ref{s2}, we will prove the Theorem \ref{the}. In section \ref{s3}, we will prove the Theorem \ref{the1}.

\section {Preliminaries}\label{s1}

In this section, we define the horizontal variables $x_{h}:=(x_{1},x_{2})$, the vertical variable $x_{v}:=x_{3}$ and $\nabla_{h}:=(\partial_{1},\partial_{2})$. For any $s,s'\in\mathbb{R}$, the anisotropic Sobolev space $H^{s,s'}$ is defined by the Sobolev space with regularity $H^{s}$ in $x_{h}$ and $H^{s'}$ in $x_{3}$. Similarly, $\dot{H}^{s,s'}$ represents  the homogenous anisotropic Sobolev space. We define $\|\cdot\|_{L^{p}}$ to denote the $L^{p}(\mathbb{R}^{3})$ norm. For exponents $p,q\in[1,\infty)$, $L_{h}^{p}(L_{v}^{q})$ denotes the space $L^{p}(\mathbb{R}_{x_{1}}\times\mathbb{R}_{x_{2}},L^{q}(\mathbb{R}_{x_{3}}))$ which is  endowed with the norm
\begin{align*}
\|u\|_{L_{h}^{p}(L_{v}^{q})}:=\{\int_{\mathbb{R}^{2}}(\int_{\mathbb{R}}|u(x_{h},x_{3})|^{q}dx_{3})^{\frac{p}{q}}dx_{h}\}^{\frac{1}{p}}.
\end{align*}
Similar, the space $L_{v}^{q}(L_{h}^{p})$ is defined.  $C$ represents genetic constant which may change from line to line. Finally, we give the following main results.

\begin{theorem}\label{the}
Let $u_{0}\in H^{0,1}(\mathbb{R}^{3})$ and $v_{0}\in H^{1}(\mathbb{R}^{3})$  and $\theta_{0}\in H^{1}(\mathbb{R}^{3})$ with $\nabla\cdot u_{0}=0$, $\alpha\geq\frac{3}{2}$ and $\beta\geq4$. The system (\ref{1}) has a unique global solution $(u,v,\theta)$ satisfying
\begin{align*}
u&\in L_{loc}^{\infty}(\mathbb{R}^{+};H^{0,1}(\mathbb{R}^{3}))\cap L_{loc}^{2}(\mathbb{R}^{+};H^{1,1}(\mathbb{R}^{3}))\cap
L_{loc}^{\beta+1}(\mathbb{R}^{+};L^{\beta+1}(\mathbb{R}^{3})),\\
\lvert u\rvert^{\frac{\beta-1}{2}}\partial_{3}u&\in L_{loc}^{2}(\mathbb{R}^{+};L^{2}(\mathbb{R}^{3})),~
\partial_{3}\lvert u\rvert^{\frac{\beta+1}{2}}\in L_{loc}^{2}(\mathbb{R}^{+};L^{2}(\mathbb{R}^{3})),\\
v&\in L_{loc}^{\infty}(\mathbb{R}^{+};H^{1}(\mathbb{R}^{3}))\cap L_{loc}^{2}(\mathbb{R}^{+};H^{1+\alpha}(\mathbb{R}^{3}))
\end{align*}
and
\begin{align*}
\theta\in L_{loc}^{\infty}(\mathbb{R}^{+};H^{1}(\mathbb{R}^{3}))\cap L_{loc}^{2}(\mathbb{R}^{+};H^{2}(\mathbb{R}^{3})).
\end{align*}
\end{theorem}

\begin{theorem}\label{the1}
Let $u_{0}\in H^{s}(\mathbb{R}^{3})$ and $v_{0}\in H^{s}(\mathbb{R}^{3})$  and $\theta_{0}\in H^{s}(\mathbb{R}^{3})$ with  $\nabla\cdot u_{0}=0$ and $s>2$, $\alpha\geq\frac{3}{2}$ and $4\leq\beta\leq5$.
The system (\ref{1}) has a unique global smooth solution $(u,v,\theta)$ satisfying
\begin{align*}
u&\in L_{loc}^{\infty}(\mathbb{R}^{+};H^{s}(\mathbb{R}^{3})),~~~\nabla_{h}u\in L_{loc}^{2}(\mathbb{R}^{+};H^{s}(\mathbb{R}^{3})),\\
v&\in L_{loc}^{\infty}(\mathbb{R}^{+};H^{s}(\mathbb{R}^{3}))\cap L_{loc}^{2}(\mathbb{R}^{+};H^{s+\alpha}(\mathbb{R}^{3}))
\end{align*}
and
\begin{align*}
\theta\in L_{loc}^{\infty}(\mathbb{R}^{+};H^{s}(\mathbb{R}^{3}))\cap L_{loc}^{2}(\mathbb{R}^{+};H^{s+1}(\mathbb{R}^{3})).
\end{align*}
\end{theorem}

\section{Proof of Theorem \ref{the}}\label{s2}
In this section, we will prove the Theorem \ref{the}.\\
Proof. We only prove the case $\frac{3}{2}\leq\alpha<\frac{5}{2}$ and $\beta\geq4$, and the
case $\alpha\geq\frac{5}{2}$ and $\beta\geq4$ can be disposed similarly. Multiplying the first equation of (\ref{1}) by $u$, the second equation of (\ref{1}) by $v$, the third equation of (\ref{1}) by $\theta$, respectively, integrating their results on $\mathbb{R}^{3}$ and integration by parts, we get
\begin{align*}
&\frac{1}{2}\frac{d}{dt}(\|u\|^{2}_{L^{2}}+\|v\|^{2}_{L^{2}}+\|\theta\|^{2}_{L^{2}})+\|\nabla_{h}u\|^{2}_{L^{2}}+\|\Lambda^{\alpha}v\|^{2}_{L^{2}}
+\|\nabla \theta\|^{2}_{L^{2}}+\|u\|^{\beta+1}_{L^{\beta+1}}\nonumber\\
&=-\int_{\mathbb{R}^{3}}\nabla\cdot(v\otimes v)\cdot udx-\int_{\mathbb{R}^{3}}(v\cdot\nabla)u\cdot vdx-\int_{\mathbb{R}^{3}}\nabla\theta\cdot vdx
-\int_{\mathbb{R}^{3}}(\nabla\cdot v)\theta dx\nonumber\\
&=0.
\end{align*}
Integrating the above equality in $[0,t]$, we have
\begin{align}\label{2}
&\|u\|^{2}_{L^{2}}+\|v\|^{2}_{L^{2}}+\|\theta\|^{2}_{L^{2}}+2\int_{0}^{t}(\|\nabla_{h}u\|^{2}_{L^{2}}+\|\Lambda^{\alpha}v\|^{2}_{L^{2}}
+\|\nabla \theta\|^{2}_{L^{2}}+\|u\|^{\beta+1}_{L^{\beta+1}})ds\nonumber\\
&=\|u_{0}\|^{2}_{L^{2}}+\|v_{0}\|^{2}_{L^{2}}+\|\theta_{0}\|^{2}_{L^{2}}.
\end{align}
Multiplying the first equation of (\ref{1}) by $-\partial_{3}^{2}u$ and integrating over on $\mathbb{R}^{3}$ and integration by parts, we have
\begin{align}\label{5}
&\frac{1}{2}\frac{d}{dt}\|\partial_{3}u\|^{2}_{L^{2}}+\|\nabla_{h}\partial_{3}u\|^{2}_{L^{2}}
+\|\lvert u\rvert^{\frac{\beta-1}{2}}\partial_{3}u\|^{2}_{L^{2}}
+\frac{4(\beta-1)}{(\beta+1)^{2}}\|\partial_{3}\lvert u\rvert^{\frac{\beta+1}{2}}\|^{2}_{L^{2}}\nonumber\\
&=\int_{\mathbb{R}^{3}}(u\cdot\nabla)u\cdot \partial_{3}^{2}udx+\int_{\mathbb{R}^{3}}\nabla\cdot(v\otimes v)\cdot\partial_{3}^{2}udx\nonumber\\
&:=I_{1}(t)+I_{2}(t).
\end{align}
For the term $I_{1}(t)$,  by integration by parts and $\nabla\cdot u=0$, it yields
\begin{align*}
I_{1}(t)&=-\sum\limits_{k,l=1}^{3}\int_{\mathbb{R}^{3}}\partial_{3}u_{k}\partial_{k}u_{l}\partial_{3}u_{l}dx\notag\\
&=-\sum\limits_{k=1}^{2}\sum\limits_{l=1}^{3}
\int_{\mathbb{R}^{3}}\partial_{3}u_{k}\partial_{k}u_{l}\partial_{3}u_{l}dx
+\sum\limits_{l=1}^{3}
\int_{\mathbb{R}^{3}}(\partial_{1}u_{1}+\partial_{2}u_{2})\partial_{3}u_{l}\partial_{3}u_{l}dx\notag\\
&:=I_{11}(t)+I_{12}(t).
\end{align*}
For the term $I_{11}(t)$, by the integration by parts, we have for $\beta>3$
\begin{align}\label{3}
\lvert I_{11}(t)\rvert&\leq \sum\limits_{k=1}^{2}\sum\limits_{l=1}^{3}
\int_{\mathbb{R}^{3}}\lvert u_{l}\rvert\lvert \partial_{3}u_{k}\rvert^{\frac{2}{\beta-1}}\lvert\partial_{3}u_{k}\rvert^{\frac{\beta-3}{\beta-1}}
\lvert\partial_{k}\partial_{3}u_{l}\rvert dx\notag\\
&+\sum\limits_{k=1}^{2}\sum\limits_{l=1}^{3}
\int_{\mathbb{R}^{3}}\lvert u_{l}\rvert\lvert\partial_{3}u_{l}\rvert^{\frac{2}{\beta-1}}\lvert\partial_{3}u_{l}\rvert^{\frac{\beta-3}{\beta-1}}
\lvert\partial_{k}\partial_{3}u_{k}\rvert dx\notag\\
&\leq \sum\limits_{k=1}^{2}\sum\limits_{l=1}^{3}\|\lvert u_{l}\rvert\lvert\partial_{3}u_{k}\rvert^{\frac{2}{\beta-1}}\|_{L^{\beta-1}}
\|\lvert\partial_{3}u_{k}\rvert^{\frac{\beta-3}{\beta-1}}\|_{L^{\frac{2(\beta-1)}{\beta-3}}}
\|\partial_{k}\partial_{3}u_{l}\|_{L^{2}}\notag\\
&+\sum\limits_{k=1}^{2}\sum\limits_{l=1}^{3}\|\lvert u_{l}\rvert\lvert\partial_{3}u_{l}\rvert^{\frac{2}{\beta-1}}\|_{L^{\beta-1}}
\|\lvert\partial_{3}u_{l}\rvert^{\frac{\beta-3}{\beta-1}}\|_{L^{\frac{2(\beta-1)}{\beta-3}}}
\|\partial_{k}\partial_{3}u_{k}\|_{L^{2}}\notag\\
&\leq \frac{1}{4}\||u|^{\frac{\beta-1}{2}}\partial_{3}u\|^{2}_{L^{2}}
+\frac{1}{4}\|\nabla_{h}\partial_{3}u\|^{2}_{L^{2}}
+C\|\partial_{3}u\|^{2}_{L^{2}}.
\end{align}
For the term $I_{12}(t)$, by H\"{o}lder and Young inequalities, we have for $\beta>3$
\begin{align}
\lvert I_{12}(t)\rvert&\leq2\sum\limits_{k=1}^{2}\sum\limits_{l=1}^{3}
\int_{\mathbb{R}^{3}}\lvert u_{k}\rvert\lvert\partial_{3}u_{l}\rvert^{\frac{2}{\beta-1}}\lvert\partial_{3}u_{l}\rvert^{\frac{\beta-3}{\beta-1}}
\lvert\partial_{k}\partial_{3}u_{l}\rvert dx\notag\\
&\leq \frac{1}{4}\|\lvert u\rvert^{\frac{\beta-1}{2}}\partial_{3}u\|^{2}_{L^{2}}
+\frac{1}{4}\|\nabla_{h}\partial_{3}u\|^{2}_{L^{2}}
+C\|\partial_{3}u\|^{2}_{L^{2}}.
\end{align}
For the term $I_{2}(t)$, using the integration by parts and H\"{o}lder inequality and Gagliardo-Nirenberg inequality and Young inequality, we deduce for
$\frac{5}{4}\leq\alpha<\frac{5}{2}$
\begin{align}\label{4}
I_{2}(t)&\leq C\|\nabla^{2}v\|_{L^{\frac{6}{5-2\alpha}}}\|v\|_{L^{\frac{3}{\alpha-1}}}\|\partial_{3}u\|_{L^{2}}
+C\|\nabla v\|_{L^{4}}^{2}\|\partial_{3}u\|_{L^{2}}\nonumber\\
&\leq C\|\Lambda^{1+\alpha}v\|_{L^{2}}\|v\|_{L^{\frac{3}{\alpha-1}}}\|\partial_{3}u\|_{L^{2}}\nonumber\\
&\leq \frac{1}{8}\|\Lambda^{1+\alpha}v\|_{L^{2}}^{2}+C\|v\|_{L^{\frac{3}{\alpha-1}}}^{2}\|\partial_{3}u\|_{L^{2}}^{2}\nonumber\\
&\leq \frac{1}{8}\|\Lambda^{1+\alpha}v\|_{L^{2}}^{2}+C\|v\|_{L^{2}}^{\frac{4\alpha-5}{\alpha}}
\|\Lambda^{\alpha}v\|_{L^{2}}^{\frac{5-2\alpha}{\alpha}}\|\partial_{3}u\|_{L^{2}}^{2}\nonumber\\
&\leq \frac{1}{8}\|\Lambda^{1+\alpha}v\|_{L^{2}}^{2}+C(\|v\|_{L^{2}}^{2}+
\|\Lambda^{\alpha}v\|_{L^{2}}^{2})\|\partial_{3}u\|_{L^{2}}^{2}.
\end{align}
Inserting (\ref{3})-(\ref{4}) into (\ref{5}), it yields for
$\frac{5}{4}\leq\alpha<\frac{5}{2}$ and  $\beta>3$
\begin{align}\label{9}
&\frac{1}{2}\frac{d}{dt}\|\partial_{3}u\|^{2}_{L^{2}}+\frac{1}{2}\|\nabla_{h}\partial_{3}u\|^{2}_{L^{2}}+\frac{1}{2}\|\lvert u\rvert^{\frac{\beta-1}{2}}\partial_{3}u\|^{2}_{L^{2}}
+\frac{4(\beta-1)}{(\beta+1)^{2}}\|\partial_{3}\lvert u\rvert^{\frac{\beta+1}{2}}\|^{2}_{L^{2}}\nonumber\\
&\leq\frac{1}{8}\|\Lambda^{1+\alpha}v\|_{L^{2}}^{2}+C(1+\|v\|_{L^{2}}^{2}+
\|\Lambda^{\alpha}v\|_{L^{2}}^{2})\|\partial_{3}u\|^{2}_{L^{2}}.
\end{align}
Multiplying the second equation of (\ref{1}) by $-\Delta v$ and integrating over on $\mathbb{R}^{3}$ and integration by parts, we deduce
\begin{align}\label{6}
&\frac{1}{2}\frac{d}{dt}\|\nabla v\|^{2}_{L^{2}}+\|\Lambda^{1+\alpha}v\|^{2}_{L^{2}}\nonumber\\
&=\int_{\mathbb{R}^{3}}(u\cdot\nabla)v\cdot \Delta vdx+\int_{\mathbb{R}^{3}}(v\cdot\nabla)u\cdot\Delta vdx+
\int_{\mathbb{R}^{3}}\nabla\theta\cdot \Delta vdx\nonumber\\
&:=I_{3}(t)+I_{4}(t)+I_{5}(t).
\end{align}
For the term $I_{3}(t)$, by virtue of the H\"{o}lder inequality and Gagliardo-Nirenberg inequality and Young inequality, we get for $\frac{5}{4}\leq\alpha<\frac{5}{2}$
\begin{align}\label{7}
I_{3}(t)&\leq \|u\|_{L^{2}}\|\nabla v\|_{L^{\frac{3}{\alpha-1}}}\|\Delta v\|_{L^{\frac{6}{5-2\alpha}}}\nonumber\\
&\leq C\|u\|_{L^{2}}\|\nabla v\|_{L^{\frac{3}{\alpha-1}}}\|\Lambda^{1+\alpha}v\|_{L^{2}}\nonumber\\
&\leq C\|u\|_{L^{2}}\|\nabla v\|_{L^{2}}^{\frac{4\alpha-5}{2\alpha}}\|\Lambda^{1+\alpha}v\|_{L^{2}}^{\frac{5}{2\alpha}}\nonumber\\
&\leq \frac{1}{8}\|\Lambda^{1+\alpha}v\|_{L^{2}}^{2}+C\|u\|_{L^{2}}^{\frac{4\alpha}{4\alpha-5}}\|\nabla v\|_{L^{2}}^{2}.
\end{align}
For the term $I_{4}(t)$, by virtue of the H\"{o}lder inequality and Sobolev embeddings and Young inequality, we have for $\frac{3}{2}\leq\alpha<\frac{5}{2}$
\begin{align}
I_{4}(t)&\leq \int_{\mathbb{R}^{3}}\lvert v\rvert \lvert\nabla_{h}u\rvert \lvert\Delta v\rvert dx+\int_{\mathbb{R}^{3}}\lvert v\rvert
\lvert\partial_{3}u\rvert \lvert\Delta v\rvert dx\nonumber\\
&\leq C\|v\|_{L^{\frac{3}{\alpha-1}}}\|\nabla_{h}u\|_{L^{2}}\|\Delta v\|_{L^{\frac{6}{5-2\alpha}}}+C\|v\|_{L^{\frac{3}{\alpha-1}}}\|\partial_{3}u\|_{L^{2}}\|\Delta v\|_{L^{\frac{6}{5-2\alpha}}}\nonumber\\
&\leq C\|v\|_{L^{\frac{3}{\alpha-1}}}\|\nabla_{h}u\|_{L^{2}}\|\Lambda^{1+\alpha} v\|_{L^{2}}+C\|v\|_{L^{\frac{3}{\alpha-1}}}\|\partial_{3}u\|_{L^{2}}\|\Lambda^{1+\alpha} v\|_{L^{2}}\nonumber\\
&\leq \frac{1}{8}\|\Lambda^{1+\alpha}v\|_{L^{2}}^{2}+C\|v\|_{L^{\frac{3}{\alpha-1}}}^{2}\|\nabla_{h}u\|_{L^{2}}^{2}
+C\|v\|_{L^{\frac{3}{\alpha-1}}}^{2}\|\partial_{3}u\|_{L^{2}}^{2}\nonumber\\
&\leq \frac{1}{8}\|\Lambda^{1+\alpha}v\|_{L^{2}}^{2}+C\|v\|_{L^{2}}^{2\alpha-3}
\|\nabla v\|_{L^{2}}^{5-2\alpha}\|\nabla_{h}u\|_{L^{2}}^{2}
+C\|v\|_{L^{2}}^{\frac{4\alpha-5}{\alpha}}\|\Lambda^{\alpha}v\|_{L^{2}}^{\frac{5-2\alpha}{\alpha}}\|\partial_{3}u\|_{L^{2}}^{2}\nonumber\\
&\leq \frac{1}{8}\|\Lambda^{1+\alpha}v\|_{L^{2}}^{2}+C\|\nabla_{h}u\|_{L^{2}}^{2}(\|v\|_{L^{2}}^{2}+\|\nabla v\|_{L^{2}}^{2})
+C(\|v\|_{L^{2}}^{2}+\|\Lambda^{\alpha}v\|_{L^{2}}^{2})\|\partial_{3}u\|_{L^{2}}^{2}.
\end{align}
For the term $I_{5}(t)$, similarly, it yields for $\alpha \geq1$
\begin{align}\label{8}
I_{5}(t)&\leq \|\nabla\theta\|_{L^{2}}\|\Delta v\|_{L^{2}}\nonumber\\
&\leq C\|\nabla\theta\|_{L^{2}}\|\nabla v\|_{L^{2}}^{\frac{\alpha-1}{\alpha}}\|\Lambda^{1+\alpha}v\|_{L^{2}}^{\frac{1}{\alpha}}\nonumber\\
&\leq \frac{1}{8}\|\Lambda^{1+\alpha}v\|_{L^{2}}^{2}+C(\|\nabla\theta\|_{L^{2}}^{2}+\|\nabla v\|_{L^{2}}^{2}).
\end{align}
Inserting (\ref{7})-(\ref{8}) into (\ref{6}), we get for $\frac{3}{2}\leq\alpha<\frac{5}{2}$
\begin{align}\label{10}
&\frac{1}{2}\frac{d}{dt}\|\nabla v\|^{2}_{L^{2}}+\frac{5}{8}\|\Lambda^{1+\alpha}v\|^{2}_{L^{2}}
\leq C(1+\|u\|_{L^{2}}^{\frac{4\alpha}{4\alpha-5}}+\|\nabla_{h}u\|_{L^{2}}^{2})\|\nabla v\|_{L^{2}}^{2}\nonumber\\
&+C(\|v\|_{L^{2}}^{2}+\|\Lambda^{\alpha}v\|_{L^{2}}^{2})\|\partial_{3}u\|_{L^{2}}^{2}+C\|\nabla_{h}u\|_{L^{2}}^{2}\|v\|_{L^{2}}^{2}
+C\|\nabla\theta\|_{L^{2}}^{2}.
\end{align}
Adding (\ref{9}) and (\ref{10}), we get for $\frac{3}{2}\leq\alpha<\frac{5}{2}$ and $\beta>3$
\begin{align}
&\frac{d}{dt}(\|\partial_{3}u\|^{2}_{L^{2}}+\|\nabla v\|^{2}_{L^{2}})+\|\nabla_{h}\partial_{3}u\|^{2}_{L^{2}}
+\|\Lambda^{1+\alpha}v\|^{2}_{L^{2}}\nonumber\\
&+\|\lvert u\rvert^{\frac{\beta-1}{2}}\partial_{3}u\|^{2}_{L^{2}}
+\|\partial_{3}\lvert u\rvert^{\frac{\beta+1}{2}}\|^{2}_{L^{2}}\nonumber\\
&\leq C(1+\|u\|_{L^{2}}^{\frac{4\alpha}{4\alpha-5}}+\|\nabla_{h}u\|_{L^{2}}^{2})\|\nabla v\|_{L^{2}}^{2}\nonumber\\
&+C(1+\|v\|_{L^{2}}^{2}+\|\Lambda^{\alpha}v\|_{L^{2}}^{2})\|\partial_{3}u\|_{L^{2}}^{2}+C\|\nabla_{h}u\|_{L^{2}}^{2}\|v\|_{L^{2}}^{2}
+C\|\nabla\theta\|_{L^{2}}^{2}.
\end{align}
Applying the Gronwall inequality and (\ref{2}), it is easy to get
\begin{align}\label{17}
\|\partial_{3}u\|^{2}_{L^{2}}+\|\nabla v\|^{2}_{L^{2}}&+\int_{0}^{t}(\|\nabla_{h}\partial_{3}u\|^{2}_{L^{2}}+\|\Lambda^{1+\alpha}v\|^{2}_{L^{2}}+\|\lvert u\rvert^{\frac{\beta-1}{2}}\partial_{3}u\|^{2}_{L^{2}}
+\|\partial_{3}\lvert u\rvert^{\frac{\beta+1}{2}}\|^{2}_{L^{2}})ds\nonumber\\
& \leq C(t,u_{0},v_{0},\theta_{0}).
\end{align}
Multiplying the thrid equation of (\ref{1}) by $-\Delta\theta$ and integrating over on $\mathbb{R}^{3}$ and integration by parts, it yields
\begin{align*}
&\frac{1}{2}\frac{d}{dt}\|\nabla \theta\|^{2}_{L^{2}}+\|\Delta\theta\|^{2}_{L^{2}}\nonumber\\
&=\int_{\mathbb{R}^{3}}(u\cdot\nabla)\theta\Delta\theta dx+\int_{\mathbb{R}^{3}}\nabla\cdot v\Delta\theta dx\nonumber\\
&\leq\|u\|_{L^{\beta+1}}\|\nabla\theta\|_{L^{\frac{2(\beta+1)}{\beta-1}}}\|\Delta\theta\|_{L^{2}}
+C\|v\|_{L^{2}}^{1-\frac{1}{\alpha}}\|\Lambda^{\alpha}v\|_{L^{2}}^{\frac{1}{\alpha}}\|\Delta\theta\|_{L^{2}}\nonumber\\
&\leq C\|u\|_{L^{\beta+1}}\|\nabla\theta\|_{L^{2}}^{\frac{\beta-2}{\beta+1}}\|\Delta\theta\|_{L^{2}}^{\frac{\beta+4}{\beta+1}}
+C\|v\|_{L^{2}}^{1-\frac{1}{\alpha}}\|\Lambda^{\alpha}v\|_{L^{2}}^{\frac{1}{\alpha}}\|\Delta\theta\|_{L^{2}}\nonumber\\
&\leq \frac{1}{2}\|\Delta\theta\|_{L^{2}}^{2}+C\|u\|_{L^{\beta+1}}^{\frac{2(\beta+1)}{\beta-2}}\|\nabla\theta\|_{L^{2}}^{2}
+C(\|v\|_{L^{2}}^{2}+\|\Lambda^{\alpha}v\|_{L^{2}}^{2})\nonumber\\
&\leq \frac{1}{2}\|\Delta\theta\|_{L^{2}}^{2}+C(1+\|u\|_{L^{\beta+1}}^{\beta+1})\|\nabla\theta\|_{L^{2}}^{2}
+C(\|v\|_{L^{2}}^{2}+\|\Lambda^{\alpha}v\|_{L^{2}}^{2}),
\end{align*}
here, we have used $\frac{2}{\beta-2}\leq1$ for $\beta\geq4$. Moreover, it is easy to get
\begin{align}
\frac{d}{dt}\|\nabla \theta\|^{2}_{L^{2}}+\|\Delta\theta\|^{2}_{L^{2}}\leq C(1+\|u\|_{L^{\beta+1}}^{\beta+1})\|\nabla\theta\|_{L^{2}}^{2}
+C(\|v\|_{L^{2}}^{2}+\|\Lambda^{\alpha}v\|_{L^{2}}^{2}).
\end{align}
Applying the Gronwall inequality and (\ref{2}), we deduce for $\beta\geq4$
\begin{align}\label{18}
\|\nabla \theta\|^{2}_{L^{2}}+\int_{0}^{t}\|\Delta\theta\|^{2}_{L^{2}}ds\leq C(t,u_{0},v_{0},\theta_{0}).
\end{align}
Finally, we will prove the uniqueness of solutions of system (\ref{1}) for $\frac{3}{2}\leq\alpha<\frac{5}{2}$ and $\beta\geq4$. Let $(\bar{u},\bar{v},\bar{\theta},\bar{p})$ and $(\tilde{u},\tilde{v},\tilde{\theta},\tilde{p})$ denote the two solutions of system (\ref{1}) with the same initial data. Let $u=\bar{u}-\tilde{u}$, $v=\bar{v}-\tilde{v}$, $\theta=\bar{\theta}-\tilde{\theta}$ and $p=\bar{p}-\tilde{p}$. We introduce the following system
 \begin{equation}\label{11}\begin{array}{l}
\left\{
\begin{array}{l}
\partial_{t}u+(\bar{u}\cdot\nabla)u+(u\cdot\nabla)\tilde{u}-\Delta_{h}u+\lvert \bar{u}\rvert^{\beta-1}\bar{u}-\lvert \tilde{u}\rvert^{\beta-1}\tilde{u}\\
+\nabla\cdot(\bar{v}\otimes v)+\nabla\cdot(v\otimes \tilde{v})+\nabla p=0,\\
\partial_{t}v+(\bar{u}\cdot\nabla)v+(u\cdot\nabla)\tilde{v}+(-\Delta)^{\alpha}v+(\bar{v}\cdot\nabla)u+(v\cdot\nabla)\tilde{u}+\nabla\theta=0,\\
\partial_{t}\theta+(\bar{u}\cdot\nabla)\theta+(u\cdot\nabla)\tilde{\theta}-\Delta \theta+\nabla\cdot v=0,\\
\nabla\cdot u=0.
\end{array}
\right.
\end{array}\end{equation}
Multiplying the first equation of (\ref{11}) by $u$, the second equation of (\ref{11}) by $v$, the third equation of (\ref{11}) by $\theta$, respectively, integrating their results on $\mathbb{R}^{3}$ and integration by parts, we have
\begin{align}\label{14}
&\frac{1}{2}\frac{d}{dt}(\|u\|^{2}_{L^{2}}+\|v\|^{2}_{L^{2}}+\|\theta\|^{2}_{L^{2}})+\|\nabla_{h}u\|^{2}_{L^{2}}+\|\Lambda^{\alpha}v\|^{2}_{L^{2}}
+\|\nabla \theta\|^{2}_{L^{2}}+\int_{\mathbb{R}^{3}}(\lvert \bar{u}\rvert^{\beta-1}\bar{u}-\lvert \tilde{u}\rvert^{\beta-1}\tilde{u})udx\nonumber\\
&=-\int_{\mathbb{R}^{3}}(u\cdot\nabla)\tilde{u}\cdot udx-\int_{\mathbb{R}^{3}}\nabla\cdot(\bar{v}\otimes v)\cdot udx-\int_{\mathbb{R}^{3}}
\nabla\cdot(v\otimes \tilde{v})\cdot udx\nonumber\\
&-\int_{\mathbb{R}^{3}}(u\cdot\nabla)\tilde{v}\cdot vdx-\int_{\mathbb{R}^{3}}(\bar{v}\cdot\nabla)u\cdot vdx-\int_{\mathbb{R}^{3}}
(v\cdot\nabla)\tilde{u}\cdot vdx-\int_{\mathbb{R}^{3}}\nabla\theta\cdot vdx\nonumber\\
&-\int_{\mathbb{R}^{3}}(u\cdot\nabla)\tilde{\theta}\theta dx-\int_{\mathbb{R}^{3}}\nabla\cdot v \theta dx\nonumber\\
&:=\sum\limits_{i=6}^{14}I_{i}(t).
\end{align}
Inspired by \cite{liu,liu1,liu2}, it yields
\begin{align}\label{15}
\int_{\mathbb{R}^{3}}(\lvert \bar{u}\rvert^{\beta-1}\bar{u}-\lvert \tilde{u}\rvert^{\beta-1}\tilde{u})udx\geq0.
\end{align}
For the terms $I_{7}(t)$, $I_{10}(t)$, $I_{12}(t)$ and $I_{14}(t)$, by integration by parts, we get
\begin{align}
I_{7}(t)+I_{10}(t)&=-\int_{\mathbb{R}^{3}}\nabla\cdot(\bar{v}\otimes v)\cdot udx
-\int_{\mathbb{R}^{3}}(\bar{v}\cdot\nabla)u\cdot vdx=0,\\
I_{12}(t)+I_{14}(t)&=-\int_{\mathbb{R}^{3}}\nabla\theta\cdot vdx-\int_{\mathbb{R}^{3}}\nabla\cdot v \theta dx=0.
\end{align}
By virtue of the Theorem 0.1 in \cite{bessaih,liu4}, we get for any $\psi\in L^{2}_{x_{3}}\cap \dot{H}^{1}_{h}$,
\begin{align}
\|\psi\|^{2}_{L_{x_{3}}^{2}(L_{h}^{4})}&\leq C\int_{\mathbb{R}}\|\psi\|_{L_{h}^{2}}\|\nabla_{h}\psi\|_{L_{h}^{2}}dx_{3}
\leq C\|\psi\|_{L^{2}}\|\nabla_{h}\psi\|_{L^{2}},\label{12}\\
\|\psi(\cdot,x_{3})\|_{L_{h}^{2}}^{2}&=\int_{-\infty}^{x_{3}}\frac{d}{dz}
(\|\psi(\cdot,z)\|_{L_{h}^{2}}^{2})dz\nonumber\\
&=2\int_{-\infty}^{x_{3}}\int_{\mathbb{R}^{2}}\psi(x_{h},z)\partial_{z}\psi(x_{h},z)dx_{h}dz\nonumber\\
&\leq C\|\psi\|_{L^{2}}\|\partial_{3}\psi\|_{L^{2}}.\label{13}
\end{align}
For the term $I_{6}(t)$, applying Sobolev embedding  and H\"{o}lder inequality and Young inequality and (\ref{12}) and (\ref{13}), we deduce
\begin{align}
I_{6}(t)&=-\sum\limits_{k=1}^{2}\sum\limits_{l=1}^{3}\int_{\mathbb{R}^{3}}u_{k}\partial_{k}\tilde{u}_{l}
u_{l}dx-\sum\limits_{l=1}^{3}\int_{\mathbb{R}^{3}}u_{3}\partial_{3}\tilde{u}_{l}
u_{l}dx\nonumber\\
&\leq \sum\limits_{k=1}^{2}\sum\limits_{l=1}^{3}\int_{\mathbb{R}}\|\partial_{k}\tilde{u}_{l}\|_{L^{2}_{h}}
\|u_{k}\|_{L^{4}_{h}}\|u_{l}\|_{L^{4}_{h}}dx_{3}
+\sum\limits_{l=1}^{3}\int_{\mathbb{R}}\|u_{3}\|_{L^{2}_{h}}\|\partial_{3}\tilde{u}_{l}\|_{L^{4}_{h}}
\|u_{l}\|_{L^{4}_{h}}dx_{3}\nonumber\\
&\leq C\|\nabla_{h}\tilde{u}\|_{L^{\infty}_{x_{3}}(L^{2}_{h})}
\|u\|_{L^{2}_{x_{3}}(L^{4}_{h})}^{2}+C\|u_{3}\|_{L^{\infty}_{x_{3}}(L^{2}_{h})}\|\partial_{3}\tilde{u}\|_{L^{2}_{x_{3}}
(L^{4}_{h})}\|u\|_{L^{2}_{x_{3}}(L^{4}_{h})}\nonumber\\
&\leq C\|\nabla_{h}\partial_{3}\tilde{u}\|_{L^{2}}^{\frac{1}{2}}\|\nabla_{h}\tilde{u}\|_{L^{2}}^{\frac{1}{2}}
\|u\|_{L^{2}}\|\nabla_{h}u\|_{L^{2}}\nonumber\\
&+C\|u_{3}\|_{L^{2}}^{\frac{1}{2}}
\|\partial_{3}u_{3}\|_{L^{2}}^{\frac{1}{2}}
\|\partial_{3}\tilde{u}\|_{L^{2}}^{\frac{1}{2}}\|\nabla_{h}\partial_{3}\tilde{u}\|_{L^{2}}^{\frac{1}{2}}
\|u\|_{L^{2}}^{\frac{1}{2}}\|\nabla_{h}u\|_{L^{2}}^{\frac{1}{2}}\nonumber\\
&\leq C\|\nabla_{h}\partial_{3}\tilde{u}\|_{L^{2}}^{\frac{1}{2}}\|\nabla_{h}\tilde{u}\|_{L^{2}}^{\frac{1}{2}}
\|u\|_{L^{2}}\|\nabla_{h}u\|_{L^{2}}\nonumber\\
&+C\|u_{3}\|_{L^{2}}^{\frac{1}{2}}
\|div_{h}u_{h}\|_{L^{2}}^{\frac{1}{2}}
\|\partial_{3}\tilde{u}\|_{L^{2}}^{\frac{1}{2}}\|\nabla_{h}\partial_{3}\tilde{u}\|_{L^{2}}^{\frac{1}{2}}
\|u\|_{L^{2}}^{\frac{1}{2}}\|\nabla_{h}u\|_{L^{2}}^{\frac{1}{2}}\nonumber\\
&\leq \frac{1}{4}\|\nabla_{h}u\|_{L^{2}}^{2}
+C(\|\nabla_{h}\tilde{u}\|_{L^{2}}^{2}+
\|\partial_{3}\tilde{u}\|_{L^{2}}^{2}+\|\nabla_{h}\partial_{3}\tilde{u}\|_{L^{2}}^{2})
\|u\|_{L^{2}}^{2}.
\end{align}
For the term $I_{8}(t)$, applying Sobolev embedding  and H\"{o}lder inequality and Young inequality, we have
\begin{align}
I_{8}(t)&\leq C\|v\|_{L^{2}_{x_{3}}(L^{4}_{h})}\|\nabla\tilde{v}\|_{L^{\infty}_{x_{3}}(L^{2}_{h})}\|u\|_{L^{2}_{x_{3}}(L^{4}_{h})}
+C\|\tilde{v}\|_{L^{\frac{3}{\alpha-1}}}\|\nabla v\|_{L^{\frac{6}{5-2\alpha}}}\|u\|_{L^{2}}\nonumber\\
&\leq C\|v\|_{L^{2}}^{\frac{1}{2}}\|\nabla v\|_{L^{2}}^{\frac{1}{2}}\|\nabla\tilde{v}\|_{L^{2}}^{\frac{1}{2}}\|\nabla^{2}\tilde{v}\|_{L^{2}}^{\frac{1}{2}}
\|u\|_{L^{2}}^{\frac{1}{2}}\|\nabla_{h}u\|_{L^{2}}^{\frac{1}{2}}+C\|\tilde{v}\|_{L^{\frac{3}{\alpha-1}}}\|\Lambda^{\alpha}v\|_{L^{2}}\|u\|_{L^{2}}
\nonumber\\
&\leq \frac{1}{8}\|\nabla_{h}u\|_{L^{2}}^{2}+ \frac{1}{16}\|\Lambda^{\alpha}v\|_{L^{2}}^{2}+C\|\nabla v\|_{L^{2}}^{2}+C\|v\|_{L^{2}}\|\nabla\tilde{v}\|_{L^{2}}\|\nabla^{2}\tilde{v}\|_{L^{2}}
\|u\|_{L^{2}}\nonumber\\
&+C\|\tilde{v}\|_{L^{\frac{3}{\alpha-1}}}^{2}\|u\|_{L^{2}}^{2}\nonumber\\
&\leq \frac{1}{8}\|\nabla_{h}u\|_{L^{2}}^{2}+\frac{1}{16}\|\Lambda^{\alpha}v\|_{L^{2}}^{2}
+C\|v\|_{L^{2}}^{\frac{2\alpha-2}{\alpha}}\|\Lambda^{\alpha}v\|_{L^{2}}^{\frac{2}{\alpha}}+C\|\tilde{v}\|_{L^{2}}^{\frac{4\alpha-5}{\alpha}}
\|\Lambda^{\alpha}\tilde{v}\|_{L^{2}}^{\frac{5-2\alpha}{\alpha}}\|u\|_{L^{2}}^{2}
\nonumber\\
&+C(\|u\|^{2}_{L^{2}}+\|v\|^{2}_{L^{2}})\|\nabla\tilde{v}\|_{L^{2}}^{\frac{2\alpha-1}{\alpha}}\|\Lambda^{1+\alpha}\tilde{v}\|_{L^{2}}^{\frac{1}{\alpha}}
\nonumber\\
&\leq \frac{1}{8}\|\nabla_{h}u\|_{L^{2}}^{2}+\frac{1}{8}\|\Lambda^{\alpha}v\|_{L^{2}}^{2}\nonumber\\
&+(1+\|\tilde{v}\|_{L^{2}}^{2}+\|\nabla\tilde{v}\|_{L^{2}}^{2}
+\|\Lambda^{\alpha}\tilde{v}\|_{L^{2}}^{2}+\|\Lambda^{1+\alpha}\tilde{v}\|_{L^{2}}^{2})(\|u\|^{2}_{L^{2}}+\|v\|^{2}_{L^{2}}).
\end{align}
For the term $I_{9}(t)$, similarly, it is easy to get
\begin{align}
I_{9}(t)&\leq C\|v\|_{L^{2}_{x_{3}}(L^{4}_{h})}\|\nabla\tilde{v}\|_{L^{\infty}_{x_{3}}(L^{2}_{h})}\|u\|_{L^{2}_{x_{3}}(L^{4}_{h})}
\nonumber\\
&\leq C\|v\|_{L^{2}}^{\frac{1}{2}}\|\nabla v\|_{L^{2}}^{\frac{1}{2}}\|\nabla\tilde{v}\|_{L^{2}}^{\frac{1}{2}}
\|\nabla^{2}\tilde{v}\|_{L^{2}}^{\frac{1}{2}}
\|u\|_{L^{2}}^{\frac{1}{2}}\|\nabla_{h}u\|_{L^{2}}^{\frac{1}{2}}\nonumber\\
&\leq \frac{1}{8}\|\nabla_{h}u\|_{L^{2}}^{2}+C\|\nabla v\|_{L^{2}}^{2}+C\|v\|_{L^{2}}\|\nabla\tilde{v}\|_{L^{2}}\|\nabla^{2}\tilde{v}\|_{L^{2}}
\|u\|_{L^{2}}\nonumber\\
&\leq \frac{1}{8}\|\nabla_{h}u\|_{L^{2}}^{2}+\frac{1}{4}\|\Lambda^{\alpha}v\|_{L^{2}}^{2}+C(1+\|\nabla\tilde{v}\|_{L^{2}}^{2}
+\|\Lambda^{1+\alpha}\tilde{v}\|_{L^{2}}^{2})(\|u\|^{2}_{L^{2}}+\|v\|^{2}_{L^{2}}).
\end{align}
For the term $I_{11}(t)$, by integration by parts and H\"{o}lder inequality and Young inequality, we get for $\frac{3}{2}\leq\alpha<\frac{5}{2}$
\begin{align}
I_{11}(t)&\leq C\|v\|_{L^{2}}\|\nabla v\|_{L^{\frac{6}{5-2\alpha}}}\|\tilde{u}\|_{L^{\frac{3}{\alpha-1}}}
\nonumber\\
&\leq C\|v\|_{L^{2}}\|\Lambda^{\alpha}v\|_{L^{2}}\|\tilde{u}\|_{L^{2}}^{\frac{2\alpha-3}{2}}\|\nabla\tilde{u}\|_{L^{2}}^{\frac{5-2\alpha}{2}}
\nonumber\\
&\leq \frac{1}{8}\|\Lambda^{\alpha}v\|_{L^{2}}^{2}
+C\|\tilde{u}\|_{L^{2}}^{2\alpha-3}\|\nabla\tilde{u}\|_{L^{2}}^{5-2\alpha}\|v\|_{L^{2}}^{2}
\nonumber\\
&\leq \frac{1}{8}\|\Lambda^{\alpha}v\|_{L^{2}}^{2}
+C(\|\tilde{u}\|_{L^{2}}^{2}+\|\nabla\tilde{u}\|_{L^{2}}^{2})\|v\|_{L^{2}}^{2}\nonumber\\
&\leq \frac{1}{8}\|\Lambda^{\alpha}v\|_{L^{2}}^{2}
+C(\|\tilde{u}\|_{L^{2}}^{2}+\|\nabla_{h}\tilde{u}\|_{L^{2}}^{2}+\|\partial_{3}\tilde{u}\|_{L^{2}}^{2})\|v\|_{L^{2}}^{2}.
\end{align}
For the term $I_{13}(t)$, by  H\"{o}lder inequality and Gagliardo-Nirenberg inequality and Young inequality, we get
\begin{align}\label{16}
I_{13}(t)&\leq C\|u\|_{L^{2}}\|\nabla \tilde{\theta}\|_{L^{3}}\|\theta\|_{L^{6}}
\nonumber\\
&\leq C\|u\|_{L^{2}}\|\nabla\theta\|_{L^{2}}\|\nabla \tilde{\theta}\|_{L^{2}}^{\frac{1}{2}}\|\Delta\tilde{\theta}\|_{L^{2}}^{\frac{1}{2}}\nonumber\\
&\leq \frac{1}{2}\|\nabla\theta\|_{L^{2}}^{2}+C(\|\nabla \tilde{\theta}\|_{L^{2}}^{2}+\|\Delta\tilde{\theta}\|_{L^{2}}^{2})\|u\|_{L^{2}}^{2}.
\end{align}
Inserting (\ref{15})-(\ref{16}) into (\ref{14}), it yields
\begin{align}
&\frac{d}{dt}(\|u\|^{2}_{L^{2}}+\|v\|^{2}_{L^{2}}+\|\theta\|^{2}_{L^{2}})+\|\nabla_{h}u\|^{2}_{L^{2}}+\|\Lambda^{\alpha}v\|^{2}_{L^{2}}
+\|\nabla \theta\|^{2}_{L^{2}}\nonumber\\
&\leq C(1+\|\tilde{u}\|_{L^{2}}^{2}+\|\nabla_{h}\tilde{u}\|_{L^{2}}^{2}+
\|\partial_{3}\tilde{u}\|_{L^{2}}^{2}+\|\nabla_{h}\partial_{3}\tilde{u}\|_{L^{2}}^{2}+\|\tilde{v}\|_{L^{2}}^{2}+\|\nabla\tilde{v}\|_{L^{2}}^{2}
+\|\Lambda^{\alpha}\tilde{v}\|_{L^{2}}^{2}\nonumber\\
&+\|\Lambda^{1+\alpha}\tilde{v}\|_{L^{2}}^{2}+\|\nabla \tilde{\theta}\|_{L^{2}}^{2}+\|\Delta\tilde{\theta}\|_{L^{2}}^{2}
)(\|u\|^{2}_{L^{2}}+\|v\|^{2}_{L^{2}}).
\end{align}
By virtue of the Gronwall inequality and (\ref{2}), (\ref{17}) and (\ref{18}), the uniqueness of solution is proved. This completes the proof of Theorem \ref{the}.

\section{Proof of Theorem \ref{the1}}\label{s3}
In this section, we will prove the Theorem \ref{the1}.\\
Proof. We only prove the case $\frac{3}{2}\leq\alpha<\frac{5}{2}$ and $4\leq\beta\leq5$, and the
case $\alpha\geq\frac{5}{2}$ and $4\leq\beta\leq5$ can be disposed similarly. The proof is structured in many
steps. \\
\textbf{Step 1} Multiplying the first equation of (\ref{1}) by $-\Delta u$ and integrating over $\mathbb{R}^{3}$ and by integration by parts, we deduce
\begin{align}\label{19}
&\frac{1}{2}\frac{d}{dt}\|\nabla u\|^{2}_{L^{2}}+\|\nabla\nabla_{h}u\|^{2}_{L^{2}}+\|\lvert u\rvert^{\frac{\beta-1}{2}}\nabla u\|^{2}_{L^{2}}
+\frac{4(\beta-1)}{(\beta+1)^{2}}\|\nabla\lvert u\rvert^{\frac{\beta+1}{2}}\|^{2}_{L^{2}}\nonumber\\
&=\int_{\mathbb{R}^{3}}(u\cdot\nabla)u\cdot \Delta udx+\int_{\mathbb{R}^{3}}\nabla\cdot(v\otimes v)\cdot\Delta udx\nonumber\\
&:=J_{1}(t)+J_{2}(t).
\end{align}
For the term $J_{1}(t)$, by integration by parts, we have
\begin{align}\label{20}
J_{1}(t)&=-\sum\limits_{k=1}^{3}\sum\limits_{i=1}^{3}\int_{\mathbb{R}^{3}}\partial_{k}u_{i}\partial_{i}u\partial_{k}udx\nonumber\\
&=\sum\limits_{k=1}^{3}\sum\limits_{i=1}^{3}\int_{\mathbb{R}^{3}}u\partial_{k}u_{i}\partial_{k}\partial_{i}udx\nonumber\\
&=\sum\limits_{k=1}^{3}\sum\limits_{i=1}^{2}\int_{\mathbb{R}^{3}}u\partial_{k}u_{i}\partial_{k}\partial_{i}udx
+\sum\limits_{k=1}^{2}\int_{\mathbb{R}^{3}}u\partial_{k}u_{3}\partial_{k}\partial_{3}udx
+\int_{\mathbb{R}^{3}}u\partial_{3}u_{3}\partial_{3}\partial_{3}udx\nonumber\\
&:=K_{11}(t)+K_{12}(t)+K_{13}(t).
\end{align}
For the term $K_{11}(t)$, applying the  H\"{o}lder inequality and Young inequality, it yields
\begin{align}
K_{11}(t)&\leq\int_{\mathbb{R}^{3}}\lvert u\rvert \lvert\nabla u\rvert \lvert\nabla\nabla_{h}u\rvert dx\nonumber\\
&\leq \|\lvert u\rvert\lvert\nabla u\rvert^{\frac{2}{\beta-1}}\|_{L^{\beta-1}}
\|\lvert\nabla u\rvert^{\frac{\beta-3}{\beta-1}}\|_{L^{\frac{2(\beta-1)}{\beta-3}}}
\|\nabla\nabla_{h}u\|_{L^{2}}\nonumber\\
&\leq \frac{1}{8}\|\lvert u\rvert^{\frac{\beta-1}{2}}\nabla u\|_{L^{2}}^{2}+\frac{1}{8}
\|\nabla\nabla_{h}u\|_{L^{2}}^{2}+C\|\nabla u\|_{L^{2}}^{2}.
\end{align}
For the term $K_{12}(t)$, similarly, we have
\begin{align}
K_{12}(t)&\leq\int_{\mathbb{R}^{3}}\lvert u\rvert \lvert\nabla u\rvert \lvert\nabla\nabla_{h}u\rvert dx\nonumber\\
&\leq \frac{1}{8}\|\lvert u\rvert^{\frac{\beta-1}{2}}\nabla u\|_{L^{2}}^{2}+\frac{1}{8}
\|\nabla\nabla_{h}u\|_{L^{2}}^{2}+C\|\nabla u\|_{L^{2}}^{2}.
\end{align}
For the term $K_{13}(t)$, by integration by parts and $\nabla\cdot u=0$ and the  H\"{o}lder inequality and Young inequality, we get
\begin{align}
K_{13}(t)&=
-\sum\limits_{i=1}^{2}\int_{\mathbb{R}^{3}}u\partial_{i}u_{i}\partial_{3}\partial_{3}udx\nonumber\\
&=\sum\limits_{i=1}^{2}\int_{\mathbb{R}^{3}}\partial_{3}u\partial_{i}u_{i}\partial_{3}udx
+\sum\limits_{i=1}^{2}\int_{\mathbb{R}^{3}}u\partial_{3}u\partial_{i}\partial_{3}u_{i}dx\nonumber\\
&\leq C\|\partial_{3}u\|_{L^{2}_{x_{3}}(L^{4}_{h})}^{2}\|\nabla_{h}u\|_{L^{\infty}_{x_{3}}(L^{2}_{h})}
+\int_{\mathbb{R}^{3}}\lvert u\rvert \lvert\nabla u\rvert \lvert \nabla\nabla_{h}u\rvert dx\nonumber\\
&\leq C\|\partial_{3}u\|_{L^{2}}\|\nabla_{h}\partial_{3}u\|_{L^{2}}^{\frac{3}{2}}\|\nabla_{h}u\|_{L^{2}}^{\frac{1}{2}}
+C\|\lvert u\rvert\lvert\nabla u\rvert^{\frac{2}{\beta-1}}\|_{L^{\beta-1}}
\|\lvert\nabla u\rvert^{\frac{\beta-3}{\beta-1}}\|_{L^{\frac{2(\beta-1)}{\beta-3}}}
\|\nabla\nabla_{h}u\|_{L^{2}}\nonumber\\
&\leq C\|\partial_{3}u\|_{L^{2}}\|\nabla\nabla_{h}u\|_{L^{2}}^{\frac{3}{2}}\|\nabla u\|_{L^{2}}^{\frac{1}{2}}
+C\|\lvert u\rvert^{\frac{\beta-1}{2}}\nabla u\|_{L^{2}}^{\frac{2}{\beta-1}}
\|\nabla u\|_{L^{2}}^{\frac{\beta-3}{\beta-1}}
\|\nabla\nabla_{h}u\|_{L^{2}}\nonumber\\
&\leq \frac{1}{4}\|\lvert u\rvert^{\frac{\beta-1}{2}}\nabla u\|_{L^{2}}^{2}+\frac{1}{4}
\|\nabla\nabla_{h}u\|_{L^{2}}^{2}+C(1+\|\partial_{3}u\|_{L^{2}}^{4})\|\nabla u\|_{L^{2}}^{2}.
\end{align}
For the term $J_{2}(t)$, by integration by parts and applying the H\"{o}lder inequality and Gagliardo-Nirenberg inequality and Young inequality, it is easy to get for $\frac{3}{2}\leq\alpha<\frac{5}{2}$
\begin{align}\label{21}
J_{2}(t)&\leq C\|\nabla u\|_{L^{2}}\|\nabla^{2}v\|_{L^{\frac{6}{5-2\alpha}}}\|v\|_{L^{\frac{3}{\alpha-1}}}
+C\|\nabla u\|_{L^{2}}\|\nabla v\|_{L^{4}}^{2}\nonumber\\
&\leq C\|\nabla u\|_{L^{2}}\|\Lambda^{1+\alpha}v\|_{L^{2}}\|v\|_{L^{\frac{3}{\alpha-1}}}
+C\|\nabla u\|_{L^{2}}\|\nabla v\|_{L^{2}}^{\frac{4\alpha-3}{2\alpha}}\|\Lambda^{1+\alpha}v\|_{L^{2}}^{\frac{3}{2\alpha}}\nonumber\\
&\leq \|v\|_{L^{\frac{3}{\alpha-1}}}^{2}+\|\nabla v\|_{L^{2}}^{\frac{4\alpha-3}{\alpha}}+C\|\Lambda^{1+\alpha}v\|_{L^{2}}^{2}\|\nabla u\|_{L^{2}}^{2}
+C\|\Lambda^{1+\alpha}v\|_{L^{2}}^{\frac{3}{\alpha}}\|\nabla u\|_{L^{2}}^{2}\nonumber\\
&\leq C(\|v\|_{L^{2}}^{2}+\|\nabla v\|_{L^{2}}^{2}+\|\nabla v\|_{L^{2}}^{\frac{4\alpha-3}{\alpha}})+C(1+\|\Lambda^{1+\alpha}v\|_{L^{2}}^{2})\|\nabla u\|_{L^{2}}^{2}.
\end{align}
Inserting (\ref{20})-(\ref{21}) into (\ref{19}), it yields
\begin{align}
&\frac{d}{dt}\|\nabla u\|^{2}_{L^{2}}+\|\nabla\nabla_{h}u\|^{2}_{L^{2}}+\|\lvert u\rvert^{\frac{\beta-1}{2}}\nabla u\|^{2}_{L^{2}}
+\|\nabla\lvert u\rvert^{\frac{\beta+1}{2}}\|^{2}_{L^{2}}\nonumber\\
&\leq C(\|v\|_{L^{2}}^{2}+\|\nabla v\|_{L^{2}}^{2}+\|\nabla v\|_{L^{2}}^{\frac{4\alpha-3}{\alpha}})+
C(1+\|\partial_{3}u\|_{L^{2}}^{4}+\|\Lambda^{1+\alpha}v\|_{L^{2}}^{2})\|\nabla u\|_{L^{2}}^{2}.
\end{align}
Applying Gronwall inequality and (\ref{2}) and (\ref{17}), we have
\begin{align}\label{22}
\|\nabla u\|^{2}_{L^{2}}+\int_{0}^{t}(\|\nabla\nabla_{h}u\|^{2}_{L^{2}}+\|\lvert u\rvert^{\frac{\beta-1}{2}}\nabla u\|^{2}_{L^{2}}
+\|\nabla\lvert u\rvert^{\frac{\beta+1}{2}}\|^{2}_{L^{2}})ds\leq C(t,u_{0},v_{0},\theta_{0}).
\end{align}
\textbf{Step 2}  Applying $\Delta$ to the third equation of (\ref{1}) and taking the $L^{2}$-inner product with $\Delta\theta$, by
integration by parts, we have for $\beta\geq4$
\begin{align*}
&\frac{1}{2}\frac{d}{dt}\|\Delta \theta\|^{2}_{L^{2}}+\|\Delta\nabla \theta\|^{2}_{L^{2}}\nonumber\\
&=-\int_{\mathbb{R}^{3}}\Delta(u\cdot\nabla\theta)\Delta \theta dx+\int_{\mathbb{R}^{3}}\Delta\nabla\cdot v\Delta \theta dx\nonumber\\
&\leq \|\nabla u\|_{L^{2}}\|\Delta\nabla \theta\|_{L^{2}}\|\nabla \theta\|_{L^{\infty}}
+\|u\|_{L^{\beta+1}}\|\nabla^{2}\theta\|_{L^{\frac{2(\beta+1)}{\beta-1}}}\|\Delta\nabla \theta\|_{L^{2}}
+\|\nabla^{2}v\|_{L^{2}}\|\Delta\nabla \theta\|_{L^{2}}\nonumber\\
&\leq C\|\nabla u\|_{L^{2}}\|\Delta\nabla \theta\|_{L^{2}}(\|\nabla \theta\|_{L^{2}}+\|\nabla \theta\|_{L^{2}}^{\frac{1}{4}}\|\Delta\nabla \theta\|_{L^{2}}^{\frac{3}{4}})
+\|u\|_{L^{\beta+1}}\|\Delta\theta\|_{L^{2}}^{\frac{\beta-2}{\beta+1}}\|\Delta\nabla \theta\|_{L^{2}}^{\frac{\beta+4}{\beta+1}}\nonumber\\
&+C\|\Delta v\|_{L^{2}}\|\Delta\nabla \theta\|_{L^{2}}\nonumber\\
&\leq \frac{1}{2}\|\Delta\nabla \theta\|_{L^{2}}^{2}+C\|\nabla u\|_{L^{2}}^{2}\|\nabla \theta\|_{L^{2}}^{2}+C\|\nabla u\|_{L^{2}}^{8}\|\nabla \theta\|_{L^{2}}^{2}+C\|u\|_{L^{\beta+1}}^{\frac{2(\beta+1)}{\beta-2}}\|\Delta\theta\|_{L^{2}}^{2}
+C\|\Delta v\|_{L^{2}}^{2}\nonumber\\
&\leq \frac{1}{2}\|\Delta\nabla \theta\|_{L^{2}}^{2}+C(1+\|u\|_{L^{\beta+1}}^{\beta+1})\|\Delta\theta\|_{L^{2}}^{2}+C\|\nabla u\|_{L^{2}}^{2}\|\nabla \theta\|_{L^{2}}^{2}\nonumber\\
&+C\|\nabla u\|_{L^{2}}^{8}\|\nabla \theta\|_{L^{2}}^{2}+C(\|\nabla v\|_{L^{2}}^{2}+\|\Lambda^{1+\alpha}v\|_{L^{2}}^{2}).
\end{align*}
Moreover, it yields
\begin{align}
\frac{d}{dt}\|\Delta \theta\|^{2}_{L^{2}}+\|\Delta\nabla \theta\|^{2}_{L^{2}}
&\leq C(1+\|u\|_{L^{\beta+1}}^{\beta+1})\|\Delta\theta\|_{L^{2}}^{2}+C\|\nabla u\|_{L^{2}}^{2}\|\nabla \theta\|_{L^{2}}^{2}\nonumber\\
&+C\|\nabla u\|_{L^{2}}^{8}\|\nabla \theta\|_{L^{2}}^{2}+C(\|\nabla v\|_{L^{2}}^{2}+\|\Lambda^{1+\alpha}v\|_{L^{2}}^{2}).
\end{align}
Applying Gronwall inequality and (\ref{2}), (\ref{17}), (\ref{18}) and (\ref{22}), we have
\begin{align}\label{35}
\|\Delta \theta\|^{2}_{L^{2}}+\int_{0}^{t}\|\Delta\nabla \theta\|^{2}_{L^{2}}ds
\leq C(t,u_{0},v_{0},\theta_{0}).
\end{align}
\textbf{Step 3} Multiplying the first equation of (\ref{1}) by $\partial_{t}u$, the second equation of (\ref{1}) by $\partial_{t}v$, respectively, and integrating over $\mathbb{R}^{3}$ and by integration by parts, we have
\begin{align*}
&\|\partial_{t}u\|^{2}_{L^{2}}+\|\partial_{t}v\|^{2}_{L^{2}}+\frac{1}{2}\frac{d}{dt}(\|\nabla_{h}u\|^{2}_{L^{2}}+\|\Lambda^{\alpha}v\|^{2}_{L^{2}})
+\frac{1}{\beta+1}\frac{d}{dt}\|u\|^{\beta+1}_{L^{\beta+1}}\nonumber\\
&=-\int_{\mathbb{R}^{3}}(u\cdot\nabla)u\cdot\partial_{t}udx-\int_{\mathbb{R}^{3}}\nabla\cdot(v\otimes v)\cdot\partial_{t}udx\nonumber\\
&-\int_{\mathbb{R}^{3}}(u\cdot\nabla)v\cdot\partial_{t}vdx-\int_{\mathbb{R}^{3}}(v\cdot\nabla)u\cdot\partial_{t}vdx
-\int_{\mathbb{R}^{3}}\nabla\theta\cdot\partial_{t}vdx\nonumber\\
&\leq \frac{1}{2}(\|\partial_{t}u\|^{2}_{L^{2}}+\|\partial_{t}v\|^{2}_{L^{2}})+C(\|\lvert u\rvert^{\frac{\beta-1}{2}}\nabla u\|_{L^{2}}^{2}+\|\nabla u\|_{L^{2}}^{2}+\|\nabla \theta\|_{L^{2}}^{2})\nonumber\\
&+C\|\nabla v\|_{L^{\frac{6}{5-2\alpha}}}^{2}\|v\|_{L^{\frac{3}{\alpha-1}}}^{2}
+C\|\nabla v\|_{L^{\frac{6}{5-2\alpha}}}^{2}\|u\|_{L^{\frac{3}{\alpha-1}}}^{2}
+C\|v\|_{L^{\infty}}^{2}\|\nabla u\|_{L^{2}}^{2}\nonumber\\
&\leq \frac{1}{2}(\|\partial_{t}u\|^{2}_{L^{2}}+\|\partial_{t}v\|^{2}_{L^{2}})+C(\|\lvert u\rvert^{\frac{\beta-1}{2}}\nabla u\|_{L^{2}}^{2}+\|\nabla u\|_{L^{2}}^{2}+\|\nabla \theta\|_{L^{2}}^{2})\nonumber\\
&+C(\|u\|_{L^{2}}^{2}+\|\nabla u\|_{L^{2}}^{2}+\|v\|_{L^{2}}^{2}+\|\nabla v\|_{L^{2}}^{2})\|\Lambda^{\alpha}v\|_{L^{2}}^{2}\nonumber\\
&+C(\|v\|_{L^{2}}^{2}+\|\Lambda^{1+\alpha}v\|_{L^{2}}^{2})\|\nabla u\|_{L^{2}}^{2}.
\end{align*}
Moreover, it yields
\begin{align}
&\|\partial_{t}u\|^{2}_{L^{2}}+\|\partial_{t}v\|^{2}_{L^{2}}+\frac{d}{dt}(\|\nabla_{h}u\|^{2}_{L^{2}}+\|\Lambda^{\alpha}v\|^{2}_{L^{2}}
+\|u\|^{\beta+1}_{L^{\beta+1}})
\nonumber\\
&\leq C(\|\lvert u\rvert^{\frac{\beta-1}{2}}\nabla u\|_{L^{2}}^{2}+\|\nabla u\|_{L^{2}}^{2}+\|\nabla \theta\|_{L^{2}}^{2})\nonumber\\
&+C(\|u\|_{L^{2}}^{2}+\|\nabla u\|_{L^{2}}^{2}+\|v\|_{L^{2}}^{2}+\|\nabla v\|_{L^{2}}^{2})\|\Lambda^{\alpha}v\|_{L^{2}}^{2}\nonumber\\
&+C(\|v\|_{L^{2}}^{2}+\|\Lambda^{1+\alpha}v\|_{L^{2}}^{2})\|\nabla u\|_{L^{2}}^{2}.
\end{align}
Applying Gronwall inequality and (\ref{2}), (\ref{17}), (\ref{18}) and (\ref{22}), we get
\begin{align}\label{27}
\|\nabla_{h}u\|^{2}_{L^{2}}+\|\Lambda^{\alpha}v\|^{2}_{L^{2}}
+\|u\|^{\beta+1}_{L^{\beta+1}}+\int_{0}^{t}(\|\partial_{s}u\|^{2}_{L^{2}}+\|\partial_{s}v\|^{2}_{L^{2}})ds\leq C(t,u_{0},v_{0},\theta_{0}).
\end{align}
\textbf{Step 4} Multiplying the third equation of (\ref{1}) by $\partial_{t}\theta$ and integrating over $\mathbb{R}^{3}$ and by integration by parts, we deduce
\begin{align*}
&\|\partial_{t}\theta\|^{2}_{L^{2}}+\frac{1}{2}\frac{d}{dt}\|\nabla \theta\|^{2}_{L^{2}}
\nonumber\\
&=-\int_{\mathbb{R}^{3}}(u\cdot\nabla)\theta\partial_{t}\theta dx-\int_{\mathbb{R}^{3}}\nabla\cdot v\partial_{t}\theta dx\nonumber\\
&\leq \frac{1}{2}\|\partial_{t}\theta\|^{2}_{L^{2}}+C(\|u\|_{L^{6}}^{2}\|\nabla\theta \|_{L^{3}}^{2}+\|\nabla v\|_{L^{2}}^{2})\nonumber\\
&\leq \frac{1}{2}\|\partial_{t}\theta\|^{2}_{L^{2}}+C(\|\nabla u\|_{L^{2}}^{2}\|\nabla\theta\|_{L^{2}}\|\Delta\theta\|_{L^{2}}+\|\nabla v\|_{L^{2}}^{2})\nonumber\\
&\leq \frac{1}{2}\|\partial_{t}\theta\|^{2}_{L^{2}}+C(\|\nabla u\|_{L^{2}}^{8}+\|\nabla v\|_{L^{2}}^{2}+\|\nabla\theta\|_{L^{2}}^{4}+\|\Delta\theta\|_{L^{2}}^{2}).
\end{align*}
Moreover, it yields
\begin{align}\label{23}
\|\partial_{t}\theta\|^{2}_{L^{2}}+\frac{d}{dt}\|\nabla \theta\|^{2}_{L^{2}}
\leq C(\|\nabla u\|_{L^{2}}^{8}+\|\nabla v\|_{L^{2}}^{2}+\|\nabla\theta\|_{L^{2}}^{4}+\|\Delta\theta\|_{L^{2}}^{2}).
\end{align}
Integrating (\ref{23}) in $[0,t]$ and by (\ref{17}) and (\ref{18}) and (\ref{22}), it is easy to get
\begin{align}\label{28}
\|\nabla \theta\|^{2}_{L^{2}}+\int_{0}^{t}\|\partial_{s}\theta\|^{2}_{L^{2}}ds
\leq C(t,u_{0},v_{0},\theta_{0}).
\end{align}
\textbf{Step 5}  Applying $\partial_{t}$ to the first equation of (\ref{1}) and $\partial_{t}$ to the second equation of (\ref{1}) and taking the $L^{2}$-inner product with $\partial_{t}u$ and  $\partial_{t}v$, respectively, by integration by parts, we get
\begin{align}\label{24}
&\frac{1}{2}\frac{d}{dt}(\|\partial_{t}u\|^{2}_{L^{2}}+\|\partial_{t}v\|^{2}_{L^{2}})+\|\nabla_{h}\partial_{t}u\|^{2}_{L^{2}}
+\|\Lambda^{\alpha}\partial_{t}v\|^{2}_{L^{2}}
+\int_{\mathbb{R}^{3}}\partial_{t}(\lvert u\rvert^{\beta-1}u)\partial_{t}udx\nonumber\\
&=-\int_{\mathbb{R}^{3}}(\partial_{t}u\cdot\nabla)u\cdot\partial_{t}udx-\int_{\mathbb{R}^{3}}\partial_{t}\nabla\cdot(v\otimes v)\cdot\partial_{t}udx
\nonumber\\
&-\int_{\mathbb{R}^{3}}(\partial_{t}u\cdot\nabla)v\cdot\partial_{t}vdx-\int_{\mathbb{R}^{3}}(\partial_{t}v\cdot\nabla)u\cdot\partial_{t}vdx\nonumber\\
&-\int_{\mathbb{R}^{3}}(v\cdot\nabla)\partial_{t}u\cdot\partial_{t}vdx-\int_{\mathbb{R}^{3}}\partial_{t}\nabla\theta\cdot\partial_{t}vdx\nonumber\\
&:=\sum\limits_{i=3}^{8}J_{i}(t).
\end{align}
Inspired by \cite{liux}, we get
\begin{align}\label{25}
\int_{\mathbb{R}^{3}}\partial_{t}(\lvert u\rvert^{\beta-1}u)\partial_{t}udx\geq0.
\end{align}
For the terms $J_{4}(t)$ and $J_{7}(t)$, by integration by parts, we get for $\frac{3}{2}\leq\alpha<\frac{5}{2}$
\begin{align}
&J_{4}(t)+J_{7}(t)\nonumber\\
&=\int_{\mathbb{R}^{3}}(\partial_{t}v\cdot\nabla)\partial_{t}u\cdot vdx\nonumber\\
&\leq \|\partial_{t}u\|_{L^{2}}\|\nabla\partial_{t}v\|_{L^{\frac{6}{5-2\alpha}}}\|v\|_{L^{\frac{3}{\alpha-1}}}
+\|\partial_{t}u\|_{L^{2}}\|\nabla v\|_{L^{\frac{6}{5-2\alpha}}}\|\partial_{t}v\|_{L^{\frac{3}{\alpha-1}}}\nonumber\\
&\leq C\|\partial_{t}u\|_{L^{2}}\|\Lambda^{\alpha}\partial_{t}v\|_{L^{2}}\|v\|_{L^{\frac{3}{\alpha-1}}}
+C\|\partial_{t}u\|_{L^{2}}\|\Lambda^{\alpha}v\|_{L^{2}}\|\partial_{t}v\|_{L^{2}}^{\frac{4\alpha-5}{2\alpha}}
\|\Lambda^{\alpha}\partial_{t}v\|_{L^{2}}^{\frac{5-2\alpha}{2\alpha}}\nonumber\\
&\leq \frac{1}{8}\|\Lambda^{\alpha}\partial_{t}v\|_{L^{2}}^{2}+C\|v\|_{L^{\frac{3}{\alpha-1}}}^{2}\|\partial_{t}u\|_{L^{2}}^{2}
+\|\partial_{t}u\|_{L^{2}}^{2}\nonumber\\
&+C\|\Lambda^{\alpha}v\|_{L^{2}}^{2}\|\partial_{t}v\|_{L^{2}}^{\frac{4\alpha-5}{\alpha}}
\|\Lambda^{\alpha}\partial_{t}v\|_{L^{2}}^{\frac{5-2\alpha}{\alpha}}\nonumber\\
&\leq \frac{1}{4}\|\Lambda^{\alpha}\partial_{t}v\|_{L^{2}}^{2}+C(1+\|v\|_{L^{\frac{3}{\alpha-1}}}^{2})\|\partial_{t}u\|_{L^{2}}^{2}
+C\|\Lambda^{\alpha}v\|_{L^{2}}^{\frac{4\alpha}{4\alpha-5}}\|\partial_{t}v\|_{L^{2}}^{2}\nonumber\\
&\leq \frac{1}{4}\|\Lambda^{\alpha}\partial_{t}v\|_{L^{2}}^{2}+C(1+\|v\|_{L^{2}}^{2}+\|\nabla v\|_{L^{2}}^{2})\|\partial_{t}u\|_{L^{2}}^{2}
+C\|\Lambda^{\alpha}v\|_{L^{2}}^{\frac{4\alpha}{4\alpha-5}}\|\partial_{t}v\|_{L^{2}}^{2}.
\end{align}
For the term $J_{3}(t)$, by straightforward computations and $\nabla\cdot u=0$, we deduce
\begin{align}
J_{3}(t)&=\sum\limits_{i=1}^{2}\sum\limits_{j=1}^{3}\int_{\mathbb{R}^{3}}\partial_{t}u_{i}\partial_{i}u_{j}\partial_{t}u_{j}dx
+\sum\limits_{j=1}^{3}\int_{\mathbb{R}^{3}}\partial_{t}u_{3}\partial_{3}u_{j}\partial_{t}u_{j}dx\nonumber\\
&\leq C\|\partial_{t}u\|_{L^{2}_{x_{3}}(L^{4}_{h})}^{2}\|\nabla_{h}u\|_{L^{\infty}_{x_{3}}(L^{2}_{h})}
+C\|\partial_{t}u\|_{L^{2}_{x_{3}}(L^{4}_{h})}\|\partial_{3}u\|_{L^{2}_{x_{3}}(L^{4}_{h})}\|\partial_{t}u_{3}\|_{L^{\infty}_{x_{3}}(L^{2}_{h})}\nonumber\\
&\leq C\|\partial_{t}u\|_{L^{2}}\|\nabla_{h}\partial_{t}u\|_{L^{2}}\|\nabla_{h}u\|_{L^{2}}^{\frac{1}{2}}\|\nabla_{h}\partial_{3}u\|_{L^{2}}^{\frac{1}{2}}
\nonumber\\
&+C\|\partial_{t}u\|_{L^{2}}\|\nabla_{h}\partial_{t}u\|_{L^{2}}^{\frac{1}{2}}\|\partial_{3}u\|_{L^{2}}^{\frac{1}{2}}
\|\nabla_{h}\partial_{3}u\|_{L^{2}}^{\frac{1}{2}}\|\partial_{t}\partial_{3}u_{3}\|_{L^{2}}^{\frac{1}{2}}\nonumber\\
&\leq \frac{1}{4}\|\nabla_{h}\partial_{t}u\|_{L^{2}}^{2}+C(\|\nabla_{h}u\|_{L^{2}}^{2}+\|\partial_{3}u\|_{L^{2}}^{2}+\|\nabla_{h}\partial_{3}u\|_{L^{2}}^{2})
\|\partial_{t}u\|_{L^{2}}^{2}.
\end{align}
For the term $J_{5}(t)$, applying the H\"{o}lder inequality and Gagliardo-Nirenberg inequality and Young inequality, we get
\begin{align}
J_{5}(t)&\leq \|\partial_{t}u\|_{L^{2}}\|\nabla v\|_{L^{\frac{6}{5-2\alpha}}}\|\partial_{t}v\|_{L^{\frac{3}{\alpha-1}}}\nonumber\\
&\leq C\|\partial_{t}u\|_{L^{2}}\|\Lambda^{\alpha}v\|_{L^{2}}\|\partial_{t}v\|_{L^{2}}^{\frac{4\alpha-5}{2\alpha}}
\|\Lambda^{\alpha}\partial_{t}v\|_{L^{2}}^{\frac{5-2\alpha}{2\alpha}}\nonumber\\
&\leq\frac{1}{8}\|\Lambda^{\alpha}\partial_{t}v\|_{L^{2}}^{2}+C\|\partial_{t}u\|_{L^{2}}^{2}
+C\|\Lambda^{\alpha}v\|_{L^{2}}^{\frac{4\alpha}{4\alpha-5}}\|\partial_{t}v\|_{L^{2}}^{2}.
\end{align}
For the term $J_{6}(t)$, applying the H\"{o}lder inequality and Young inequality and by integration by parts, we have
\begin{align}
J_{6}(t)&\leq C\|\partial_{t}v\|_{L^{2}}\|\nabla\partial_{t}v\|_{L^{\frac{6}{5-2\alpha}}}\|u\|_{L^{\frac{3}{\alpha-1}}}\nonumber\\
&\leq C\|\partial_{t}v\|_{L^{2}}\|\Lambda^{\alpha}\partial_{t}v\|_{L^{2}}\|u\|_{L^{\frac{3}{\alpha-1}}}\nonumber\\
&\leq\frac{1}{16}\|\Lambda^{\alpha}\partial_{t}v\|_{L^{2}}^{2}
+C(\|u\|_{L^{2}}^{2}+\|\nabla u\|_{L^{2}}^{2})\|\partial_{t}v\|_{L^{2}}^{2}.
\end{align}
For the term $J_{8}(t)$, similarly, one has
\begin{align}\label{26}
J_{8}(t)&\leq C\|\nabla\partial_{t}v\|_{L^{2}}\|\partial_{t}\theta\|_{L^{2}}\nonumber\\
&\leq\frac{1}{16}\|\Lambda^{\alpha}\partial_{t}v\|_{L^{2}}^{2}
+C(\|\partial_{t}v\|_{L^{2}}^{2}+\|\partial_{t}\theta\|_{L^{2}}^{2}).
\end{align}
Inserting (\ref{25})-(\ref{26}) into (\ref{24}), we get
\begin{align}
&\frac{d}{dt}(\|\partial_{t}u\|^{2}_{L^{2}}+\|\partial_{t}v\|^{2}_{L^{2}})+\|\nabla_{h}\partial_{t}u\|^{2}_{L^{2}}
+\|\Lambda^{\alpha}\partial_{t}v\|^{2}_{L^{2}}
\nonumber\\
&\leq C(1+\|\nabla_{h}u\|_{L^{2}}^{2}+\|\partial_{3}u\|_{L^{2}}^{2}+\|\nabla_{h}\partial_{3}u\|_{L^{2}}^{2}+\|v\|_{L^{2}}^{2}+\|\nabla v\|_{L^{2}}^{2})\|\partial_{t}u\|_{L^{2}}^{2}\nonumber\\
&+C(1+\|u\|_{L^{2}}^{2}+\|\nabla u\|_{L^{2}}^{2}+\|\Lambda^{\alpha}v\|_{L^{2}}^{\frac{4\alpha}{4\alpha-5}})\|\partial_{t}v\|_{L^{2}}^{2}
+C\|\partial_{t}\theta\|_{L^{2}}^{2}.
\end{align}
Applying Gronwall inequality and (\ref{2}), (\ref{17}), (\ref{22}), (\ref{27}) and (\ref{28}), we have
\begin{align}\label{29}
\|\partial_{t}u\|^{2}_{L^{2}}+\|\partial_{t}v\|^{2}_{L^{2}}+\int_{0}^{t}(\|\nabla_{h}\partial_{s}u\|^{2}_{L^{2}}
+\|\Lambda^{\alpha}\partial_{s}v\|^{2}_{L^{2}})ds
\leq C(t,u_{0},v_{0},\theta_{0}).
\end{align}
\textbf{Step 6}  Taking inner product of the first equation of (\ref{1}) with $-\Delta u$, it yields
\begin{align}
&\|\nabla\nabla_{h}u\|^{2}_{L^{2}}+\|\lvert u\rvert^{\frac{\beta-1}{2}}\nabla u\|^{2}_{L^{2}}
+\|\nabla\lvert u\rvert^{\frac{\beta+1}{2}}\|^{2}_{L^{2}}\nonumber\\
&\leq C(\|\partial_{t}u\|^{2}_{L^{2}}+\|(u\cdot\nabla)u\|^{2}_{L^{2}}+\|\nabla\cdot(v\otimes v)\|^{2}_{L^{2}})\nonumber\\
&\leq \frac{1}{2}\|\lvert u\rvert^{\frac{\beta-1}{2}}\nabla u\|^{2}_{L^{2}}+C(\|\partial_{t}u\|^{2}_{L^{2}}+\|\nabla u\|^{2}_{L^{2}}+\|\nabla v\|_{L^{\frac{6}{5-2\alpha}}}^{2}\|v\|^{2}_{L^{\frac{3}{\alpha-1}}})\nonumber\\
&\leq \frac{1}{2}\|\lvert u\rvert^{\frac{\beta-1}{2}}\nabla u\|^{2}_{L^{2}}+C(\|\partial_{t}u\|^{2}_{L^{2}}+\|\nabla u\|^{2}_{L^{2}}+\|\Lambda^{\alpha} v\|_{L^{2}}^{2}(\|v\|^{2}_{L^{2}}+\|\Lambda^{\alpha} v\|_{L^{2}}^{2})).
\end{align}
Moreover, by virtue of (\ref{2}), (\ref{22}), (\ref{27}) and (\ref{29}), it yields
\begin{align}\label{36}
&\|\nabla\nabla_{h}u\|^{2}_{L^{2}}+\|\lvert u\rvert^{\frac{\beta-1}{2}}\nabla u\|^{2}_{L^{2}}
+\|\nabla\lvert u\rvert^{\frac{\beta+1}{2}}\|^{2}_{L^{2}}\nonumber\\
&\leq C(\|\partial_{t}u\|^{2}_{L^{2}}+\|\nabla u\|^{2}_{L^{2}}+\|v\|^{4}_{L^{2}}+\|\Lambda^{\alpha} v\|_{L^{2}}^{4})\nonumber\\
&\leq C(t,u_{0},v_{0},\theta_{0}).
\end{align}
Since $\dot{H}^{1}\hookrightarrow L^{6}$, one has
\begin{align}\label{37}
\|u\|^{\beta+1}_{L^{3(\beta+1)}}
\leq C\|\nabla\lvert u\rvert^{\frac{\beta+1}{2}}\|^{2}_{L^{2}}
\leq C(t,u_{0},v_{0},\theta_{0}).
\end{align}
Taking inner product of the second equation of (\ref{1}) with $-\Delta v$, it yields
\begin{align}
\|\Lambda^{1+\alpha}v\|^{2}_{L^{2}}
&\leq C(\|\partial_{t}v\|^{2}_{L^{2}}+\|(u\cdot\nabla)v\|^{2}_{L^{2}}+\|(v\cdot\nabla)u\|_{L^{2}}^{2}+\|\nabla\theta\|^{2}_{L^{2}})\nonumber\\
&\leq C(\|\partial_{t}v\|^{2}_{L^{2}}+\|u\|^{2}_{L^{6}}\|\nabla v\|^{2}_{L^{3}}+\|v\|_{L^{\infty}}^{2}\|\nabla u\|_{L^{2}}^{2}+\|\nabla\theta\|^{2}_{L^{2}})\nonumber\\
&\leq C(\|\partial_{t}v\|^{2}_{L^{2}}+\|\nabla u\|^{2}_{L^{2}}\|v\|^{\frac{2\alpha-3}{\alpha}}_{L^{2}}\|\Lambda^{\alpha}v\|^{\frac{3}{\alpha}}_{L^{2}}+\|v\|_{L^{2}}^{2}\|\nabla u\|_{L^{2}}^{2}\nonumber\\
&+\|v\|_{L^{2}}^{\frac{2\alpha-1}{\alpha+1}}\|\Lambda^{1+\alpha}u\|^{\frac{3}{1+\alpha}}_{L^{2}}\|\nabla u\|_{L^{2}}^{2}+\|\nabla\theta\|^{2}_{L^{2}})
\nonumber\\
&\leq\frac{1}{2}\|\Lambda^{1+\alpha}v\|^{2}_{L^{2}}+C(\|\partial_{t}v\|^{2}_{L^{2}}+\|\nabla u\|^{2}_{L^{2}}\|v\|^{\frac{2\alpha-3}{\alpha}}_{L^{2}}\|\Lambda^{\alpha}v\|^{\frac{3}{\alpha}}_{L^{2}}+\|v\|_{L^{2}}^{2}\|\nabla u\|_{L^{2}}^{2}\nonumber\\
&+\|v\|_{L^{2}}^{2}\|\nabla u\|_{L^{2}}^{\frac{4(\alpha+1)}{2\alpha-1}}+\|\nabla\theta\|^{2}_{L^{2}}).
\end{align}
Moreover, by virtue of (\ref{2}), (\ref{18}), (\ref{22}), (\ref{27}) and (\ref{29}), it yields
\begin{align}
\|\Lambda^{1+\alpha}v\|^{2}_{L^{2}}
&\leq C(\|\partial_{t}v\|^{2}_{L^{2}}+\|\nabla u\|^{2}_{L^{2}}\|v\|^{\frac{2\alpha-3}{\alpha}}_{L^{2}}\|\Lambda^{\alpha}v\|^{\frac{3}{\alpha}}_{L^{2}}+\|v\|_{L^{2}}^{2}\|\nabla u\|_{L^{2}}^{2}\nonumber\\
&+\|v\|_{L^{2}}^{2}\|\nabla u\|_{L^{2}}^{\frac{4(\alpha+1)}{2\alpha-1}}+\|\nabla\theta\|^{2}_{L^{2}})\nonumber\\
&\leq C(t,u_{0},v_{0},\theta_{0}).
\end{align}
\textbf{Step 7} Applying $\Delta$ to the first equation of (\ref{1}) and taking the $L^{2}$-inner product with $\Delta u$, by
integration by parts, we have
\begin{align}\label{30}
&\frac{1}{2}\frac{d}{dt}\|\Delta u\|^{2}_{L^{2}}+\|\Delta\nabla_{h} u\|^{2}_{L^{2}}+\|\lvert u\rvert^{\frac{\beta-1}{2}}\Delta u\|^{2}_{L^{2}}\nonumber\\
&\leq C\int_{\mathbb{R}^{3}}\lvert u\rvert^{\beta-2}\lvert \nabla u\rvert \lvert \nabla u\rvert \lvert \Delta u\rvert dx-\int_{\mathbb{R}^{3}}\Delta((u\cdot\nabla)u)\cdot\Delta udx
-\int_{\mathbb{R}^{3}}\Delta\nabla\cdot(v\otimes v)\cdot\Delta udx.
\end{align}
Applying $\Delta$ to the second equation of (\ref{1}) and taking the $L^{2}$-inner product with $\Delta v$, by
integration by parts, we also have
\begin{align}\label{31}
&\frac{1}{2}\frac{d}{dt}\|\Delta v\|^{2}_{L^{2}}+\|\Lambda^{2+\alpha}v\|^{2}_{L^{2}}\nonumber\\
&=-\int_{\mathbb{R}^{3}}\Delta((u\cdot\nabla)v)\cdot\Delta vdx-\int_{\mathbb{R}^{3}}\Delta((v\cdot\nabla)u)\cdot\Delta vdx
-\int_{\mathbb{R}^{3}}\Delta\nabla\theta\cdot\Delta vdx.
\end{align}
Adding (\ref{30}) and (\ref{31}), it yields
\begin{align}\label{32}
&\frac{1}{2}\frac{d}{dt}(\|\Delta u\|^{2}_{L^{2}}+\|\Delta v\|^{2}_{L^{2}})+\|\Delta\nabla_{h} u\|^{2}_{L^{2}}+\|\lvert u\rvert^{\frac{\beta-1}{2}}\Delta u\|^{2}_{L^{2}}+\|\Lambda^{2+\alpha}v\|^{2}_{L^{2}}\nonumber\\
&\leq C\int_{\mathbb{R}^{3}}\lvert u\rvert^{\beta-2}\lvert \nabla u\rvert \lvert \nabla u\rvert \lvert \Delta u\rvert dx-\int_{\mathbb{R}^{3}}\Delta((u\cdot\nabla)u)\cdot\Delta udx
-\int_{\mathbb{R}^{3}}\Delta\nabla\cdot(v\otimes v)\cdot\Delta udx\nonumber\\
&-\int_{\mathbb{R}^{3}}\Delta((u\cdot\nabla)v)\cdot\Delta vdx-\int_{\mathbb{R}^{3}}\Delta((v\cdot\nabla)u)\cdot\Delta vdx
-\int_{\mathbb{R}^{3}}\Delta\nabla\theta\cdot\Delta vdx\nonumber\\
&:=\sum\limits_{i=9}^{14}J_{i}(t).
\end{align}
For the term $J_{9}(t)$, applying the H\"{o}lder inequality and Young inequality, we have
\begin{align}\label{33}
&J_{9}(t)\nonumber\\
&\leq C\|\lvert u\rvert^{\beta-2}\nabla u\|_{L^{2}}\|\Delta u\|_{L^{2}_{x_{3}}(L^{4}_{h})}\|\nabla u\|_{L^{\infty}_{x_{3}}(L^{4}_{h})}\nonumber\\
&\leq C\|u\|_{L^{3(\beta+1)}}^{\beta-2}\|\nabla u\|_{L^{\frac{6(\beta+1)}{\beta+7}}}\|\Delta u\|_{L^{2}}^{\frac{1}{2}}\|\Delta\nabla_{h}u\|_{L^{2}}^{\frac{1}{2}}\|\nabla u\|_{L^{2}}^{\frac{1}{4}}\|\nabla \partial_{3}u\|_{L^{2}}^{\frac{1}{4}}\|\nabla\nabla_{h}u\|_{L^{2}}^{\frac{1}{4}}\|\nabla\nabla_{h}\partial_{3}u\|_{L^{2}}^{\frac{1}{4}}\nonumber\\
&\leq C\|u\|_{L^{3(\beta+1)}}^{\beta-2}\|\Delta u\|_{L^{2}}^{\frac{3}{4}+\frac{\beta-2}{\beta+1}}\|\Delta\nabla_{h}u\|_{L^{2}}^{\frac{3}{4}}\|\nabla u\|_{L^{2}}^{\frac{1}{4}+\frac{3}{\beta+1}}\|\nabla\nabla_{h}u\|_{L^{2}}^{\frac{1}{4}}\nonumber\\
&\leq \frac{1}{8}\|\Delta\nabla_{h}u\|_{L^{2}}^{2}+C\|u\|_{L^{3(\beta+1)}}^{\frac{8(\beta-2)}{5}}\|\nabla u\|_{L^{2}}^{\frac{2(\beta+13)}{5(\beta+1)}}\|\nabla\nabla_{h}u\|_{L^{2}}^{\frac{2}{5}}\|\Delta u\|_{L^{2}}^{\frac{6}{5}+\frac{8(\beta-2)}{5(\beta+1)}}\nonumber\\
&\leq \frac{1}{8}\|\Delta\nabla_{h}u\|_{L^{2}}^{2}+C\|u\|_{L^{3(\beta+1)}}^{\frac{8(\beta-2)}{5}}\|\nabla u\|_{L^{2}}^{\frac{2(\beta+13)}{5(\beta+1)}}\|\nabla\nabla_{h}u\|_{L^{2}}^{\frac{2}{5}}(1+\|\Delta u\|_{L^{2}}^{2}),
\end{align}
here, we have used $\frac{\beta-2}{\beta+1}\leq\frac{1}{2}$ for $\beta\leq5$. For the term $J_{10}(t)$, by straightforward computations and integration by parts and $\nabla\cdot u=0$, we get
\begin{align*}
J_{10}(t)&=\sum\limits_{k=1}^{2}\int_{\mathbb{R}^{3}}\partial_{k}u\nabla u\partial_{k}\Delta udx
+\sum\limits_{k=1}^{2}\int_{\mathbb{R}^{3}}u\nabla\partial_{k}u\partial_{k}\Delta udx
-\int_{\mathbb{R}^{3}}\nabla\partial_{3}(u\nabla u)\nabla\partial_{3}udx\nonumber\\
&=\sum\limits_{k=1}^{2}\int_{\mathbb{R}^{3}}\partial_{k}u\nabla u\partial_{k}\Delta udx
+\sum\limits_{k=1}^{2}\int_{\mathbb{R}^{3}}u\nabla\partial_{k}u\partial_{k}\Delta udx
-\int_{\mathbb{R}^{3}}\nabla\partial_{3}u\nabla u\nabla\partial_{3}udx\nonumber\\
&-\int_{\mathbb{R}^{3}}\partial_{3}u\nabla\nabla u\nabla\partial_{3}udx
-\int_{\mathbb{R}^{3}}\nabla u\nabla\partial_{3}u\nabla\partial_{3}udx
\nonumber\\
&=K_{21}(t)+K_{22}(t)+K_{23}(t)+K_{24}(t)+K_{25}(t).
\end{align*}
For the terms $K_{21}(t)$ and $K_{22}(t)$, applying the H\"{o}lder inequality and Young inequality, we get
\begin{align}
K_{21}(t)+K_{22}(t)&\leq \|\nabla_{h}u\|_{L^{3}}\|\nabla u\|_{L^{6}}\|\Delta\nabla_{h}u\|_{L^{2}}+
\|u\|_{L^{6}}\|\nabla\nabla_{h}u\|_{L^{3}}\|\Delta\nabla_{h}u\|_{L^{2}}\nonumber\\
&\leq C\|\nabla_{h}u\|_{L^{2}}^{\frac{1}{2}}\|\nabla\nabla_{h}u\|_{L^{2}}^{\frac{1}{2}}\|\Delta u\|_{L^{2}}\|\Delta\nabla_{h}u\|_{L^{2}}\nonumber\\
&+
C\|\nabla u\|_{L^{2}}\|\nabla\nabla_{h}u\|_{L^{2}}^{\frac{1}{2}}\|\Delta\nabla_{h}u\|_{L^{2}}^{\frac{3}{2}}\nonumber\\
&\leq C\|\nabla_{h}u\|_{L^{2}}^{\frac{1}{2}}\|\nabla\nabla_{h}u\|_{L^{2}}^{\frac{1}{2}}\|\Delta u\|_{L^{2}}\|\Delta\nabla_{h}u\|_{L^{2}}\nonumber\\
&+
C\|\nabla u\|_{L^{2}}\|\Delta u\|_{L^{2}}^{\frac{1}{2}}\|\Delta\nabla_{h}u\|_{L^{2}}^{\frac{3}{2}}\nonumber\\
&\leq \frac{1}{8}\|\Delta\nabla_{h}u\|_{L^{2}}^{2}+C(\|\nabla_{h}u\|_{L^{2}}^{2}+\|\nabla u\|_{L^{2}}^{4}+\|\nabla\nabla_{h}u\|_{L^{2}}^{2})\|\Delta u\|_{L^{2}}^{2}.
\end{align}
For the term $K_{23}(t)$,  by straightforward computations and $\nabla\cdot u=0$ and Young inequality, we have
\begin{align}
K_{23}(t)&=-\sum\limits_{i=1}^{2}\int_{\mathbb{R}^{3}}\nabla\partial_{3}u_{i}\partial_{i}u\nabla\partial_{3}udx
-\int_{\mathbb{R}^{3}}\nabla\partial_{3}u_{3}\partial_{3}u\nabla\partial_{3}udx
\nonumber\\
&\leq C\|\nabla_{h}u\|_{L^{\infty}_{x_{3}}(L^{2}_{h})}\|\nabla\partial_{3}u\|_{L^{2}_{x_{3}}(L^{4}_{h})}^{2}
+C\|\nabla\nabla_{h}u\|_{L^{\infty}_{x_{3}}(L^{2}_{h})}\|\partial_{3}u\|_{L^{2}_{x_{3}}(L^{4}_{h})}\|\nabla\partial_{3}u\|_{L^{2}_{x_{3}}(L^{4}_{h})}\nonumber\\
&\leq C\|\Delta u\|_{L^{2}}\|\Delta \nabla_{h}u\|_{L^{2}}\|\nabla_{h}u\|_{L^{2}}^{\frac{1}{2}}\|\nabla_{h}\partial_{3}u\|_{L^{2}}^{\frac{1}{2}}\nonumber\\
&+C\|\Delta u\|_{L^{2}}\|\Delta \nabla_{h}u\|_{L^{2}}\|\partial_{3}u\|_{L^{2}}^{\frac{1}{2}}\|\nabla_{h}\partial_{3}u\|_{L^{2}}^{\frac{1}{2}}\nonumber\\
&\leq\frac{1}{16}\|\Delta \nabla_{h}u\|_{L^{2}}^{2}+C(\|\nabla_{h}u\|_{L^{2}}^{2}+\|\partial_{3}u\|_{L^{2}}^{2}+\|\nabla_{h}\partial_{3}u\|_{L^{2}}^{2})\|\Delta u\|_{L^{2}}^{2}.
\end{align}
For the term $K_{24}(t)$, similarly, it yields
\begin{align}
K_{24}(t)&=-\sum\limits_{i=1}^{2}\int_{\mathbb{R}^{3}}\partial_{3}u_{i}\nabla\partial_{i}u\nabla\partial_{3}udx
-\int_{\mathbb{R}^{3}}\partial_{3}u_{3}\nabla\partial_{3}u\nabla\partial_{3}udx
\nonumber\\
&\leq C\|\nabla\nabla_{h}u\|_{L^{\infty}_{x_{3}}(L^{2}_{h})}\|\partial_{3}u\|_{L^{2}_{x_{3}}(L^{4}_{h})}\|\nabla\partial_{3}u\|_{L^{2}_{x_{3}}(L^{4}_{h})}
+C\|\nabla_{h}u\|_{L^{\infty}_{x_{3}}(L^{2}_{h})}\|\nabla\partial_{3}u\|_{L^{2}_{x_{3}}(L^{4}_{h})}^{2}\nonumber\\
&\leq C\|\Delta u\|_{L^{2}}\|\Delta \nabla_{h}u\|_{L^{2}}\|\partial_{3}u\|_{L^{2}}^{\frac{1}{2}}\|\nabla_{h}\partial_{3}u\|_{L^{2}}^{\frac{1}{2}}\nonumber\\
&+C\|\Delta u\|_{L^{2}}\|\Delta \nabla_{h}u\|_{L^{2}}\|\nabla_{h}u\|_{L^{2}}^{\frac{1}{2}}\|\nabla_{h}\partial_{3}u\|_{L^{2}}^{\frac{1}{2}}
\nonumber\\
&\leq\frac{1}{16}\|\Delta \nabla_{h}u\|_{L^{2}}^{2}+C(\|\nabla_{h}u\|_{L^{2}}^{2}+\|\partial_{3}u\|_{L^{2}}^{2}+\|\nabla_{h}\partial_{3}u\|_{L^{2}}^{2})\|\Delta u\|_{L^{2}}^{2}.
\end{align}
For the last term $K_{25}(t)$, similarly, we also have
\begin{align}
K_{25}(t)&=-\sum\limits_{i=1}^{2}\int_{\mathbb{R}^{3}}\nabla u_{i}\partial_{i}\partial_{3}u\nabla\partial_{3}udx
-\int_{\mathbb{R}^{3}}\nabla u_{3}\partial_{3}\partial_{3}u\nabla\partial_{3}udx
\nonumber\\
&\leq C\|\nabla u\|_{L^{2}_{x_{3}}(L^{4}_{h})}\|\nabla_{h}\partial_{3}u\|_{L^{\infty}_{x_{3}}(L^{2}_{h})}\|\nabla\partial_{3}u\|_{L^{2}_{x_{3}}(L^{4}_{h})}\nonumber\\
&+C\|\nabla u_{3}\|_{L^{\infty}_{x_{3}}(L^{2}_{h})}\|\partial_{3}\partial_{3}u\|_{L^{2}_{x_{3}}(L^{4}_{h})}\|\nabla\partial_{3}u\|_{L^{2}_{x_{3}}(L^{4}_{h})}\nonumber\\
&\leq C\|\Delta u\|_{L^{2}}\|\Delta \nabla_{h}u\|_{L^{2}}\|\nabla u\|_{L^{2}}^{\frac{1}{2}}\|\nabla\nabla_{h}u\|_{L^{2}}^{\frac{1}{2}}\nonumber\\
&+C\|\Delta u\|_{L^{2}}\|\Delta \nabla_{h}u\|_{L^{2}}\|\nabla u\|_{L^{2}}^{\frac{1}{2}}\|\nabla_{h}\partial_{3}u_{3}\|_{L^{2}}^{\frac{1}{2}}
\nonumber\\
&\leq\frac{1}{8}\|\Delta \nabla_{h}u\|_{L^{2}}^{2}+C(\|\nabla u\|_{L^{2}}^{2}+\|\nabla\nabla_{h}u\|_{L^{2}}^{2})\|\Delta u\|_{L^{2}}^{2}.
\end{align}
For the term $J_{11}(t)$, applying the H\"{o}lder inequality and Young inequality, we get
\begin{align}
J_{11}(t)
&\leq C\|\Delta\nabla v\|_{L^{\frac{6}{5-2\alpha}}}\|v\|_{L^{\frac{3}{\alpha-1}}}\|\Delta u\|_{L^{2}}
+C\|\nabla^{2}v\|_{L^{\frac{6}{5-2\alpha}}}\|\nabla v\|_{L^{\frac{3}{\alpha-1}}}\|\Delta u\|_{L^{2}}\nonumber\\
&\leq C\|v\|_{L^{\frac{3}{\alpha-1}}}\|\Lambda^{2+\alpha}v\|_{L^{2}}\|\Delta u\|_{L^{2}}
+C\|\nabla v\|_{L^{\frac{3}{\alpha-1}}}\|\Lambda^{1+\alpha}v\|_{L^{2}}\|\Delta u\|_{L^{2}}\nonumber\\
&\leq C\|v\|_{L^{2}}^{\frac{4\alpha-5}{2\alpha}}\|\Lambda^{\alpha}v\|_{L^{2}}^{\frac{5-2\alpha}{2\alpha}}\|\Lambda^{2+\alpha}v\|_{L^{2}}\|\Delta u\|_{L^{2}}\nonumber\\
&+C\|\nabla v\|_{L^{2}}^{\frac{4\alpha-5}{2\alpha}}\|\Lambda^{1+\alpha}v\|_{L^{2}}^{\frac{5-2\alpha}{2\alpha}}\|\Lambda^{1+\alpha}v\|_{L^{2}}\|\Delta u\|_{L^{2}}\nonumber\\
&\leq \frac{1}{8}\|\Lambda^{2+\alpha}v\|_{L^{2}}^{2}+C(\|v\|_{L^{2}}^{2}+\|\Lambda^{\alpha}v\|_{L^{2}}^{2}+\|\Lambda^{1+\alpha}v\|_{L^{2}}^{2})\|\Delta u\|_{L^{2}}^{2}\nonumber\\
&+C(\|\nabla v\|_{L^{2}}^{2}+\|\Lambda^{1+\alpha}v\|_{L^{2}}^{2}).
\end{align}
For the term $J_{12}(t)$, by integration by parts and the H\"{o}lder inequality and Young inequality, we get
\begin{align}
J_{12}(t)
&\leq C\|\Delta\nabla v\|_{L^{\frac{6}{5-2\alpha}}}\|\nabla v\|_{L^{\frac{3}{\alpha-1}}}\|\nabla u\|_{L^{2}}
+C\|\Delta\nabla v\|_{L^{\frac{6}{5-2\alpha}}}\|\Delta v\|_{L^{2}}\|u\|_{L^{\frac{3}{\alpha-1}}}\nonumber\\
&\leq C\|\Lambda^{2+\alpha}v\|_{L^{2}}\|\nabla v\|_{L^{\frac{3}{\alpha-1}}}\|\nabla u\|_{L^{2}}
+C\|\Lambda^{2+\alpha}v\|_{L^{2}}\|\Delta v\|_{L^{2}}\|u\|_{L^{\frac{3}{\alpha-1}}}\nonumber\\
&\leq C\|\Lambda^{2+\alpha}v\|_{L^{2}}\|\nabla u\|_{L^{2}}\|\nabla v\|_{L^{2}}^{\frac{4\alpha-5}{2\alpha}}\|\Lambda^{1+\alpha} v\|_{L^{2}}^{\frac{5-2\alpha}{2\alpha}}\nonumber\\
&+C\|\Lambda^{2+\alpha}v\|_{L^{2}}\|\Delta v\|_{L^{2}}\|u\|_{L^{2}}^{\frac{2\alpha-3}{2}}\|\nabla u\|_{L^{2}}^{\frac{5-2\alpha}{2}}\nonumber\\
&\leq \frac{1}{8}\|\Lambda^{2+\alpha}v\|_{L^{2}}^{2}+C\|\Lambda^{1+\alpha} v\|_{L^{2}}^{2}+C\|\nabla v\|_{L^{2}}^{2}\|\nabla u\|_{L^{2}}^{\frac{4\alpha}{4\alpha-5}}\nonumber\\
&+C(\|u\|_{L^{2}}^{2}+\|\nabla u\|_{L^{2}}^{2})\|\Delta v\|_{L^{2}}^{2}.
\end{align}
For the term $J_{13}(t)$, similarly, it yields
\begin{align}
J_{13}(t)
&\leq C\|\Delta\nabla v\|_{L^{\frac{6}{5-2\alpha}}}\|\nabla v\|_{L^{\frac{3}{\alpha-1}}}\|\nabla u\|_{L^{2}}
+C\|\Delta\nabla v\|_{L^{\frac{6}{5-2\alpha}}}\|\Delta u\|_{L^{2}}\|v\|_{L^{\frac{3}{\alpha-1}}}\nonumber\\
&\leq C\|\Lambda^{2+\alpha}v\|_{L^{2}}\|\nabla v\|_{L^{\frac{3}{\alpha-1}}}\|\nabla u\|_{L^{2}}
+C\|\Lambda^{2+\alpha}v\|_{L^{2}}\|\Delta u\|_{L^{2}}\|v\|_{L^{\frac{3}{\alpha-1}}}\nonumber\\
&\leq C\|\Lambda^{2+\alpha}v\|_{L^{2}}\|\nabla u\|_{L^{2}}\|\nabla v\|_{L^{2}}^{\frac{4\alpha-5}{2\alpha}}\|\Lambda^{1+\alpha} v\|_{L^{2}}^{\frac{5-2\alpha}{2\alpha}}\nonumber\\
&+C\|\Lambda^{2+\alpha}v\|_{L^{2}}\|\Delta u\|_{L^{2}}\|v\|_{L^{2}}^{\frac{2\alpha-3}{2}}\|\nabla v\|_{L^{2}}^{\frac{5-2\alpha}{2}}\nonumber\\
&\leq \frac{1}{4}\|\Lambda^{2+\alpha}v\|_{L^{2}}^{2}+C\|\Lambda^{1+\alpha} v\|_{L^{2}}^{2}+C\|\nabla v\|_{L^{2}}^{2}\|\nabla u\|_{L^{2}}^{\frac{4\alpha}{4\alpha-5}}\nonumber\\
&+C(\|v\|_{L^{2}}^{2}+\|\nabla v\|_{L^{2}}^{2})\|\Delta u\|_{L^{2}}^{2}.
\end{align}
For the term $J_{14}(t)$, similarly, we also have
\begin{align}\label{34}
J_{14}(t)
&\leq C\|\Delta\nabla \theta\|_{L^{2}}\|\Delta v\|_{L^{2}}\nonumber\\
&\leq C(\|\Delta\nabla \theta\|_{L^{2}}^{2}+\|\Delta v\|_{L^{2}}^{2}).
\end{align}
Inserting (\ref{33})-(\ref{34}) into (\ref{32}), we get
\begin{align}
&\frac{d}{dt}(\|\Delta u\|^{2}_{L^{2}}+\|\Delta v\|^{2}_{L^{2}})+\|\Delta\nabla_{h} u\|^{2}_{L^{2}}+\|\lvert u\rvert^{\frac{\beta-1}{2}}\Delta u\|^{2}_{L^{2}}+\|\Lambda^{2+\alpha}v\|^{2}_{L^{2}}\nonumber\\
&\leq C\|u\|_{L^{3(\beta+1)}}^{\frac{8(\beta-2)}{5}}\|\nabla u\|_{L^{2}}^{\frac{2(\beta+13)}{5(\beta+1)}}\|\nabla\nabla_{h}u\|_{L^{2}}^{\frac{2}{5}}(1+\|\Delta u\|_{L^{2}}^{2})+C(\|\nabla u\|_{L^{2}}^{2}+\|\nabla u\|_{L^{2}}^{4}\nonumber\\
&+\|\nabla_{h}\partial_{3}u\|_{L^{2}}^{2}+\|\nabla\nabla_{h}u\|_{L^{2}}^{2}
+\|v\|_{L^{2}}^{2}+\|\nabla v\|_{L^{2}}^{2}+\|\Lambda^{\alpha}v\|_{L^{2}}^{2}+\|\Lambda^{1+\alpha}v\|_{L^{2}}^{2})\|\Delta u\|_{L^{2}}^{2}\nonumber\\
&+C(\|\nabla v\|_{L^{2}}^{2}+\|\Lambda^{1+\alpha}v\|_{L^{2}}^{2})+C\|\nabla v\|_{L^{2}}^{2}\|\nabla u\|_{L^{2}}^{\frac{4\alpha}{4\alpha-5}}\nonumber\\
&+C(1+\|u\|_{L^{2}}^{2}+\|\nabla u\|_{L^{2}}^{2})\|\Delta v\|_{L^{2}}^{2}+C\|\Delta\nabla \theta\|_{L^{2}}^{2}.
\end{align}
Applying the Gronwall inequality and (\ref{2}), (\ref{17}), (\ref{22}), (\ref{35}), (\ref{36}) and (\ref{37}), we get
\begin{align}\label{41}
\|\Delta u\|^{2}_{L^{2}}+\|\Delta v\|^{2}_{L^{2}}+\int_{0}^{t}(\|\Delta\nabla_{h} u\|^{2}_{L^{2}}+\|\Lambda^{2+\alpha}v\|^{2}_{L^{2}})ds\leq C(t,u_{0},v_{0},\theta_{0}).
\end{align}
Moreover, it yields
\begin{align}\label{42}
\|u\|_{L^{\infty}}^{2}\leq C(\|u\|^{2}_{L^{2}}+\|\Delta u\|^{2}_{L^{2}})\leq C(t,u_{0},v_{0},\theta_{0}).
\end{align}
\textbf{Step 8} Let $s>2$. Applying $\Lambda^{s}$ to the first equation of (\ref{1}), the second equation of (\ref{1}) and the third equation of (\ref{1}) and taking the $L^{2}$-inner product with $\Lambda^{s}u$, $\Lambda^{s}v$ and $\Lambda^{s}\theta$, respectively,  by integration by parts, we have
\begin{align}\label{38}
&\frac{1}{2}\frac{d}{dt}(\|\Lambda^{s}u\|^{2}_{L^{2}}+\|\Lambda^{s}v\|^{2}_{L^{2}}+\|\Lambda^{s}\theta\|^{2}_{L^{2}})
+\|\Lambda^{s}\nabla_{h}u\|^{2}_{L^{2}}+\|\Lambda^{s+\alpha}v\|^{2}_{L^{2}}+\|\Lambda^{s+1}\theta\|^{2}_{L^{2}}\nonumber\\
&=\int_{\mathbb{R}^{3}}\Lambda^{s}(u\cdot\nabla u)\cdot\Lambda^{s}udx-\int_{\mathbb{R}^{3}}\Lambda^{s}(\lvert u\rvert^{\beta-1}u)\cdot\Lambda^{s}udx
-\int_{\mathbb{R}^{3}}\Lambda^{s}\nabla\cdot(v\otimes v)\cdot\Lambda^{s}udx\nonumber\\
&-\int_{\mathbb{R}^{3}}\Lambda^{s}(u\cdot\nabla v)\cdot\Lambda^{s}vdx-\int_{\mathbb{R}^{3}}\Lambda^{s}(v\cdot\nabla u)\cdot\Lambda^{s}vdx-\int_{\mathbb{R}^{3}}\Lambda^{s}\nabla\theta\cdot\Lambda^{s}vdx
\nonumber\\
&-\int_{\mathbb{R}^{3}}\Lambda^{s}(u\cdot\nabla\theta)\Lambda^{s}\theta dx-\int_{\mathbb{R}^{3}}\Lambda^{s}\nabla\cdot v\Lambda^{s}\theta dx\nonumber\\
&:=\sum\limits_{i=15}^{22}J_{i}(t).
\end{align}
For the term $J_{15}(t)$, applying the H\"{o}lder inequality and Young inequality, we have
\begin{align}\label{39}
J_{15}(t)&\leq C\|[\Lambda^{s},u\cdot\nabla]u\|_{L^{2}}\|\Lambda^{s}u\|_{L^{2}}\nonumber\\
&\leq C\|\Lambda^{s}u\|_{L^{2}_{x_{3}}(L^{4}_{h})}\|\nabla u\|_{L^{\infty}_{x_{3}}(L^{4}_{h})}\|\Lambda^{s}u\|_{L^{2}}\nonumber\\
&\leq C\|\Lambda^{s}u\|_{L^{2}}^{\frac{3}{2}}\|\Lambda^{s}\nabla_{h}u\|_{L^{2}}^{\frac{1}{2}}\|\nabla u\|_{L^{2}}^{\frac{1}{4}}\|\nabla\nabla_{h} u\|_{L^{2}}^{\frac{1}{4}}\|\nabla \partial_{3}u\|_{L^{2}}^{\frac{1}{4}}\|\nabla\nabla_{h}\partial_{3}u\|_{L^{2}}^{\frac{1}{4}}\nonumber\\
&\leq C\|\Lambda^{s}u\|_{L^{2}}^{\frac{3}{2}}\|\Lambda^{s}\nabla_{h}u\|_{L^{2}}^{\frac{1}{2}}\|\nabla u\|_{L^{2}}^{\frac{1}{4}}\|\nabla\nabla_{h} u\|_{L^{2}}^{\frac{1}{4}}\|\Delta u\|_{L^{2}}^{\frac{1}{4}}\|\Delta\nabla_{h}u\|_{L^{2}}^{\frac{1}{4}}\nonumber\\
&\leq \frac{1}{2}\|\Lambda^{s}\nabla_{h}u\|_{L^{2}}^{2}+C(\|\nabla u\|_{L^{2}}+\|\nabla\nabla_{h} u\|_{L^{2}}+\|\Delta u\|_{L^{2}}^{2}+\|\Delta\nabla_{h}u\|_{L^{2}}^{2})\|\Lambda^{s}u\|_{L^{2}}^{2}.
\end{align}
For the term $J_{16}(t)$, similarly, we have
\begin{align}
J_{16}(t)&\leq C\|u\|_{L^{\infty}}^{\beta-1}\|\Lambda^{s}u\|_{L^{2}}^{2}.
\end{align}
For the term $J_{17}(t)$, by  Sobolev embeddings and Young inequality, we get
\begin{align}
J_{17}(t)&\leq C\|\Lambda^{s}u\|_{L^{2}}\|\Lambda^{s}\nabla v\|_{L^{\frac{6}{5-2\alpha}}}\|v\|_{L^{\frac{3}{\alpha-1}}}
+C\|\Lambda^{s}u\|_{L^{2}}\|\Lambda^{s}v\|_{L^{\frac{3}{\alpha-1}}}\|\nabla v\|_{L^{\frac{6}{5-2\alpha}}}\nonumber\\
&\leq C\|\Lambda^{s}u\|_{L^{2}}\|\Lambda^{s+\alpha}v\|_{L^{2}}\|v\|_{L^{\frac{3}{\alpha-1}}}
+C\|\Lambda^{s}u\|_{L^{2}}\|\Lambda^{s}v\|_{L^{2}}^{\frac{4\alpha-5}{2\alpha}}\|\Lambda^{s+\alpha}v\|_{L^{2}}^{\frac{5-2\alpha}{2\alpha}}
\|\Lambda^{\alpha}v\|_{L^{2}}\nonumber\\
&\leq\frac{1}{8}\|\Lambda^{s+\alpha}v\|_{L^{2}}^{2}+C(1+\|v\|_{L^{2}}^{2}+\|\Lambda^{\alpha}v\|_{L^{2}}^{2})\|\Lambda^{s}u\|_{L^{2}}^{2}
+C\|\Lambda^{\alpha}v\|_{L^{2}}^{\frac{4\alpha}{4\alpha-5}}\|\Lambda^{s}v\|_{L^{2}}^{2}.
\end{align}
For the term $J_{18}(t)$, similarly, it yields
\begin{align}
J_{18}(t)&\leq C\|\Lambda^{s}v\|_{L^{2}}\|\Lambda^{s}\nabla v\|_{L^{\frac{6}{5-2\alpha}}}\|u\|_{L^{\frac{3}{\alpha-1}}}
+C\|\Lambda^{s}u\|_{L^{2}}\|\Lambda^{s}v\|_{L^{\frac{3}{\alpha-1}}}\|\nabla v\|_{L^{\frac{6}{5-2\alpha}}}\nonumber\\
&\leq C\|\Lambda^{s}v\|_{L^{2}}\|\Lambda^{s+\alpha}v\|_{L^{2}}\|u\|_{L^{\frac{3}{\alpha-1}}}
+C\|\Lambda^{s}u\|_{L^{2}}\|\Lambda^{s}v\|_{L^{2}}^{\frac{4\alpha-5}{2\alpha}}\|\Lambda^{s+\alpha}v\|_{L^{2}}^{\frac{5-2\alpha}{2\alpha}}
\|\Lambda^{\alpha}v\|_{L^{2}}\nonumber\\
&\leq\frac{1}{8}\|\Lambda^{s+\alpha}v\|_{L^{2}}^{2}+C(\|u\|_{L^{2}}^{2}+\|\nabla u\|_{L^{2}}^{2})\|\Lambda^{s}v\|_{L^{2}}^{2}+C\|\Lambda^{s}u\|_{L^{2}}^{2}\nonumber\\
&+C\|\Lambda^{\alpha}v\|_{L^{2}}^{\frac{4\alpha}{4\alpha-5}}\|\Lambda^{s}v\|_{L^{2}}^{2}.
\end{align}
For the term $J_{19}(t)$, by straightforward computations and similar method, we have
\begin{align}
J_{19}(t)&=
-\int_{\mathbb{R}^{3}}\Lambda^{s}\nabla\cdot(v\otimes u)\cdot\Lambda^{s}vdx+\int_{\mathbb{R}^{3}}\Lambda^{s}(\nabla\cdot v\cdot u)\cdot\Lambda^{s}vdx\nonumber\\
&=\int_{\mathbb{R}^{3}}\Lambda^{s}(v\otimes u)\cdot\Lambda^{s}\nabla vdx+\int_{\mathbb{R}^{3}}\Lambda^{s}(\nabla\cdot v\cdot u)\cdot\Lambda^{s}vdx\nonumber\\
&\leq C\|\Lambda^{s}v\|_{L^{2}}\|\Lambda^{s}\nabla v\|_{L^{\frac{6}{5-2\alpha}}}\|u\|_{L^{\frac{3}{\alpha-1}}}
+C\|\Lambda^{s}u\|_{L^{2}}\|\Lambda^{s}\nabla v\|_{L^{\frac{6}{5-2\alpha}}}\|v\|_{L^{\frac{3}{\alpha-1}}}\nonumber\\
&+C\|\Lambda^{s}u\|_{L^{2}}\|\Lambda^{s}v\|_{L^{\frac{3}{\alpha-1}}}\|\nabla v\|_{L^{\frac{6}{5-2\alpha}}}\nonumber\\
&\leq C\|\Lambda^{s}v\|_{L^{2}}\|\Lambda^{s+\alpha}v\|_{L^{2}}\|u\|_{L^{\frac{3}{\alpha-1}}}
+C\|\Lambda^{s}u\|_{L^{2}}\|\Lambda^{s}v\|_{L^{2}}^{\frac{4\alpha-5}{2\alpha}}\|\Lambda^{s+\alpha}v\|_{L^{2}}^{\frac{5-2\alpha}{2\alpha}}
\|\Lambda^{\alpha}v\|_{L^{2}}\nonumber\\
&+C\|\Lambda^{s}u\|_{L^{2}}\|\Lambda^{s+\alpha}v\|_{L^{2}}\|v\|_{L^{\frac{3}{\alpha-1}}}\nonumber\\
&\leq\frac{1}{4}\|\Lambda^{s+\alpha}v\|_{L^{2}}^{2}+C(\|u\|_{L^{2}}^{2}+\|\nabla u\|_{L^{2}}^{2})\|\Lambda^{s}v\|_{L^{2}}^{2}+C(1+\|v\|_{L^{2}}^{2}+\|\nabla v\|_{L^{2}}^{2})\|\Lambda^{s}u\|_{L^{2}}^{2}\nonumber\\
&+C\|\Lambda^{\alpha}v\|_{L^{2}}^{\frac{4\alpha}{4\alpha-5}}\|\Lambda^{s}v\|_{L^{2}}^{2}.
\end{align}
For the terms $J_{20}(t)$ and $J_{22}(t)$, by integration by parts, we have
\begin{align}
J_{20}(t)+J_{22}(t)=-\int_{\mathbb{R}^{3}}\Lambda^{s}\nabla\theta\cdot\Lambda^{s}vdx+\int_{\mathbb{R}^{3}}\Lambda^{s}v\cdot\Lambda^{s}\nabla\theta dx
=0.
\end{align}
For the term $J_{21}(t)$, similarly, it yields
\begin{align}\label{40}
J_{21}(t)&\leq C\|u\|_{L^{\infty}}\|\Lambda^{s}\nabla\theta\|_{L^{2}}\|\Lambda^{s}\theta\|_{L^{2}}
+C\|\Lambda^{s}u\|_{L^{2}}\|\nabla\theta\|_{L^{3}}\|\Lambda^{s}\theta\|_{L^{6}}\nonumber\\
&\leq C\|u\|_{L^{\infty}}\|\Lambda^{s}\nabla\theta\|_{L^{2}}\|\Lambda^{s}\theta\|_{L^{2}}
+C\|\Lambda^{s}u\|_{L^{2}}\|\nabla\theta\|_{L^{2}}^{\frac{1}{2}}\|\Delta\theta\|_{L^{2}}^{\frac{1}{2}}
\|\Lambda^{s}\nabla\theta\|_{L^{2}}\nonumber\\
&\leq \frac{1}{2}\|\Lambda^{s}\nabla\theta\|_{L^{2}}^{2}+C\|u\|_{L^{\infty}}^{2}\|\Lambda^{s}\theta\|_{L^{2}}^{2}
+C(\|\nabla\theta\|_{L^{2}}^{2}+\|\Delta\theta\|_{L^{2}}^{2})\|\Lambda^{s}u\|_{L^{2}}^{2}.
\end{align}
Inserting (\ref{39})-(\ref{40}) into (\ref{38}), it is easy to get
\begin{align}
&\frac{d}{dt}(\|\Lambda^{s}u\|^{2}_{L^{2}}+\|\Lambda^{s}v\|^{2}_{L^{2}}+\|\Lambda^{s}\theta\|^{2}_{L^{2}})
+\|\Lambda^{s}\nabla_{h}u\|^{2}_{L^{2}}+\|\Lambda^{s+\alpha}v\|^{2}_{L^{2}}+\|\Lambda^{s+1}\theta\|^{2}_{L^{2}}\nonumber\\
&\leq C(1+\|u\|_{L^{2}}^{2}+\|u\|_{L^{\infty}}^{\beta-1}+\|\nabla u\|_{L^{2}}^{2}+\|\nabla\nabla_{h} u\|_{L^{2}}^{2}+\|\Delta u\|_{L^{2}}^{2}+\|\Delta\nabla_{h}u\|_{L^{2}}^{2}
\nonumber\\
&+\|v\|_{L^{2}}^{2}+\|\nabla v\|_{L^{2}}^{2}+\|\Lambda^{\alpha}v\|_{L^{2}}^{2}+\|\Lambda^{\alpha}v\|_{L^{2}}^{\frac{4\alpha}{4\alpha-5}}
+\|\nabla\theta\|_{L^{2}}^{2}+\|\Delta\theta\|_{L^{2}}^{2})\nonumber\\
&(\|\Lambda^{s}u\|_{L^{2}}^{2}+\|\Lambda^{s}v\|_{L^{2}}^{2}+\|\Lambda^{s}\theta\|_{L^{2}}^{2}).
\end{align}
Applying Gronwall inequality and (\ref{2}), (\ref{17}), (\ref{18}), (\ref{22}), (\ref{27}), (\ref{36}), (\ref{41}) and (\ref{42}), we get
\begin{align}
\|\Lambda^{s}u\|^{2}_{L^{2}}+\|\Lambda^{s}v\|^{2}_{L^{2}}+\|\Lambda^{s}\theta\|^{2}_{L^{2}}
&+\int_{0}^{t}(\|\Lambda^{s}\nabla_{h}u\|^{2}_{L^{2}}+\|\Lambda^{s+\alpha}v\|^{2}_{L^{2}}+\|\Lambda^{s+1}\theta\|^{2}_{L^{2}})\nonumber\\
&\leq C(t,u_{0},v_{0},\theta_{0}).
\end{align}
This completes the proof of Theorem \ref{the1}.
\begin{remark}
By using the similar method for the uniqueness  of Theorem \ref{the}, we prove the uniqueness  of Theorem \ref{the1}.
\end{remark}

\end{document}